\documentclass[preprint]{elsarticle}
\usepackage[a4paper, top=2cm, bottom=2cm, left=2.5cm, right=2.5cm]{geometry}
\usepackage{bm}
\usepackage{amssymb}
\usepackage{amsmath}
\usepackage{float}
\usepackage{graphicx}    
\usepackage{subcaption}  
\usepackage{booktabs}    
\usepackage{caption}     
\usepackage{multirow}    
\usepackage{array}       
\usepackage{makecell}    
\usepackage{nomencl} 
\usepackage[numbers]{natbib}

\makenomenclature
\journal{Elsevier}

\begin{document}

\begin{frontmatter}



\title{Evolutionary de-homogenization using a generative model for optimizing solid-porous infill structures considering the stress concentration issue}


\author{Shuzhi Xu, Hiroki Kawabe, Kentaro Yaji*} 

\affiliation{organization={Department of Mechanical Engineering, Graduate School of Engineering, Osaka University},
            addressline={2-1, Yamadaoka}, 
            city={Suita},
            postcode={565-0871}, 
            state={Osaka},
            country={Japan}}

\begin{abstract}
The design of porous infill structures presents significant challenges due to their complex geometric configurations, such as the accurate representation of geometric boundaries and the control of localized maximum stress. In current mainstream design methods, such as topology optimization, the analysis is often performed using pixel or voxel-based element approximations. These approximations, constrained by the optimization framework, result in substantial geometric discrepancies between the analysis model and the final physical model. Such discrepancies can severely impact structural performance, particularly for localized properties like stress response, where accurate geometry is critical to mitigating stress concentration. To address these challenges, we propose \textit{evolutionary de-homogenization}, which is a design framework based on the integration of de-homogenization and data-driven multifidelity optimization. This framework facilitates the hybrid solid-porous infill design by bridging the gap between low-fidelity analysis and high-fidelity physical realizations, ensuring both geometric accuracy and enhanced structural performance. The low-fidelity level utilizes commonly used density control variables, while the high-fidelity level involves stress analysis based on structures with precise geometric representations. By employing a de-homogenization-based mapping method, a side-by-side correspondence between low-fidelity and high-fidelity results is established. The low-fidelity control variables are iteratively adjusted to optimize the high-fidelity results by integrating deep generative model with multi-objective evolutionary algorithm. Finally, numerical experiments demonstrate the effectiveness of the proposed method.

\end{abstract}



\begin{keyword}
porous infill design\sep 
de-homogenization\sep 
stress concentration\sep 
data-driven\sep 
generative model\sep 
multifidelity optimization



\end{keyword}

\end{frontmatter}


\section{Introduction} \label{sec:1}
Porous structures, due to their superior mechanical properties such as ultra-lightweight, high stiffness, and energy absorption capabilities, have garnered widespread attention in practical applications \cite{sun2017topological}\cite{zhang2015bioinspired}. Concurrently, advancements in Additive Manufacturing (AM) have significantly accelerated the development of porous structures by reducing material usage and shortening production time \cite{guo2021recent}\cite{liu2018current}\cite{brackett2011topology}\cite{pyka2012surface}. In recent years, topology optimization \cite{bendsoe2013topology}\cite{eschenauer2001topology}\cite{sigmund2013topology}, a powerful computational tool for determining the optimal distribution of materials within a design domain, has been widely applied in the efficient design of porous structures. Despite these advances, purely porous infill designs often exhibit inferior mechanical performance compared to solid designs under the same conditions. Therefore, when emphasizing the advantages of porous infill designs, the focus often shifts toward their scalability for multifunctional applications or their potential for further weight reduction in fixed geometric configurations. To balance weight reduction and mechanical performance, the hybrid solid-porous design approach has emerged. This method improves both local and overall performance by introducing solid components into critical regions of the lattice structure. Researchers have implemented and validated these designs using topology optimization techniques, and results show that their performance surpasses that of purely porous infill designs~\cite{chen2021topology}\cite{bai2023topology}\cite{zhang2023multiscale}\cite{gao2024multi}.

Currently, there are several methods for generating porous infill structures. A classical method is multi-scale topology optimization, which is closely integrated with the traditional topology optimization framework \cite{bendsoe1988generating}\cite{wu2021topology}. In multi-scale optimization, the analysis model is divided into macro and micro scales. Each element in the macro model represents a microstructure, which is modeled as a heterogeneous material unit using homogenization theory~\cite{xia2015design}\cite{andreassen2014determine}\cite{dong2019149}. Many studies have explored how to optimize structures on both scales simultaneously~\cite{sivapuram2016simultaneous}\cite{gao2019concurrent}. However, this method has encountered bottlenecks due to its high computational cost, leading to a slowdown in research. To address this issue, various compromise methods have emerged, such as parameterized unit cells with single or multiple parameters~\cite{cheng2019functionally}\cite{wang2020concurrent}\cite{wang2018concurrent} and multi-domain parallel optimization~\cite{gao2019topology}\cite{li2018topology}\cite{li2019new}\cite{xu2021multi}. Recently, with the development of deep learning models, multi-scale optimization has regained attention under the concept of free material design~\cite{wang2023deep}\cite{wang2020deep}\cite{chu2024exploring}\cite{jia2024modulate}.

Compared to traditional solid structure optimization only in macroscale, multi-scale optimization provides greater design freedom, allowing the approach to reach theoretical limits and enabling non-conventional designs. However, this increased design freedom often results in overly complex geometries that are difficult to manufacture using mainstream methods. For instance, structures generated through multi-scale design often contain many microstructures with extremely small geometric features, making precise manufacturing more complex \cite{wu2021topology}. Additionally, due to their inconsistent multi-scale nature, these complex structures are difficult to interpret as coherent and manufacturable designs \cite{wu2021topology}. A simple solution is to scale each microstructure to the size of bilinear elements and apply appropriate orientation. The process of generating coherent and physically feasible designs through homogenization-based optimization is called de-homogenization. 
As a representative de-homogenization technique, Pantz and Trabelsi \cite{pantz2008post} proposed a post-processig method that provides an implicit geometric description for obtaining well-connected single-scale designs from spatially varying multi-scale designs. Recently, this method has been simplified and improved~\cite{groen2018homogenization}\cite{woldseth2024phasor}, and it has been applied to various problems~\cite{li2024analytical}\cite{li2024multi}\cite{feppon2024multiscale}.

Another method involves imposing local volume constraints to control local material distribution. Wu et al.~\cite{wu2017infill} introduced upper bounds on the proximity of each design element, with constraints on localized material volume fraction, and aggregated all local per-voxel constraints using the P-norm. This idea has been extended to implement maximum length scale constraints in topology optimization, achieving bone-like porous structures. The application of localized material volume constraints in generating porous infill designs has since been widely discussed and applied to various scenarios, such as self-supporting infill structures \cite{liu2021topology}, multi-material structures \cite{li2020spatial}, stress-based designs \cite{wei2024isogeometric}, and buckling designs \cite{liu2023topology}.

To fully explore the weight reduction potential of sandwich structures while ensuring mechanical performance, shell-infill features have been incorporated into the topology design model. To the best of the authors' knowledge, Clausen et al. \cite{clausen2015topology} were the first to combine the density method with a topology design method for sandwich structures, utilizing a two-step density filter to describe the shell and infill domains. Subsequently, Luo et al. \cite{luo2019topology} proposed an erosion-based interface identification method to distinguish between the shell and infill layers in sandwich structures. Yoon and Yi \cite{yoon2019new} later developed a new density filter, which generates the sandwich structure by multiplying the modified density design variables with the original design variables. Besides the density method, such designs can also be realized through other approaches \cite{wang2018level}\cite{hoang2020topology} and can be combined with the aforementioned porous design methods to achieve more complex designs \cite{groen2019homogenization}\cite{zhou2022concurrent}\cite{wu2017minimum}\cite{hoang2020adaptive}\cite{li2024topology}.

Significant challenges remain with current infill design methods using mainstream approaches. While traditional multi-scale optimization provides the greatest design freedom and the most rigorous theoretical foundation, the full-scale structures are often overly complex, making geometric representation difficult, especially in freeform design domains. In contrast, methods based on local volume constraints can accurately generate porous structures in macroscale. However, as a full-scale design method, they also require higher computational cost when representing complex geometries. Moreover, ensuring that the macroscopic structure adheres to local volume constraints often results in trade-offs with the target performance of the structure. De-homogenization provides a compromise solution for these methods, but as a post-processing technique, it inevitably introduces discrepancies between the final full-scale design and the homogenized results, particularly in terms of local performance.

We believe that an optimization framework incorporating multifidelity formulation and non-gradient optimization methods may effectively address these issues. Specifically, the multifidelity approach can quickly explore the design space using low-fidelity models, significantly reducing the number of high-fidelity model evaluations. The non-gradient optimization method does not rely on continuous geometric or gradient information, thereby avoiding complex sensitivity calculations. The combination of both approaches makes it easier to find a better solution in free-form design spaces. Recently, Yaji et al. \cite{yaji2022data} proposed a multifidelity topology design (MFTD) method in conjunction with genetic algorithm (GA)-based topology optimization. In MFTD, topology optimization problems that are difficult to solve analytically or numerically can be addressed indirectly through simplified surrogate optimization problems, while maintaining accurate performance evaluation throughout the process. Additionally, the optimization of control variables is achieved by integrating multi-objective optimization algorithms with deep generative models. MFTD has demonstrated its effectiveness across a broad range of applications, including stress minimization problem, heat exchangers, and heat energy storage design~\cite{kato2024maximum}\cite{kobayashi2019freeform}\cite{yaji2019framework}\cite{LUO2025124596}. The method has also seen several improvements, such as enhanced crossover strategies~\cite{yaji2024latent}, extension to three-dimensional cases~\cite{yang2024data}, and the introduction of multiple design variable fields~\cite{kawabe2024data}.

This paper proposes a new structural design framework, named as \textit{evolutionary de-homogenization}, based on the MFTD method for generating hybrid solid-porous infill structures to address several key challenges in porous infill design: 1. The issue of local stress concentration in porous infill designs; 2. The difficulty in accurately optimizing smooth boundaries in complex porous infill stress problems during the optimization process; 3. The significant geometric and performance gap between the final optimized results (such as the 0-1 density field in density-based methods) and the actual expressed results (the smoothed \textit{STereoLithography} (STL) model) in traditional methods. Inspired by multifidelity optimization, we first extract simplified low-fidelity data from complex solid-porous infill structures and propose a design method that maps this low-fidelity data to the final complex Computer-Aided Design (CAD) model (i.e., the detail full-scale sturcture used for high-fidelity simulation). Based on this, we optimize the low-fidelity data using an evolutionary algorithm (EA), where the performance metrics for each design variable are determined by an adaptive mesh model generated from the corresponding detail CAD model. The current optimal results are selected using the Non-Dominated Sorting Genetic Algorithm II (NSGA-II). To address the challenges of crossover operations in the EA within this framework, we employ a Variational AutoEncoder (VAE)~\cite{kingma2013auto}. After multiple iterations of selection, the optimized CAD model is ultimately obtained, which can be accurately converted into an STL model for further AM.

The remainder of this paper is organized as follows: Section~\ref*{sec:Multifidelity Method} introduces the design method for the solid-porous infill structure, detailing how a simplified low-fidelity design variable is abstracted from the complex design for use in MFTD. It also provides a comprehensive explanation of the multifidelity formulation, including the low-fidelity optimization and high-fidelity evaluation, with a focus on the mapping process from low-fidelity data to high-fidelity results. Section~\ref*{sec:3} presents the flowchart of the proposed multifidelity optimization framework, illustrating how the optimization of low-fidelity variables leads to the optimization of high-fidelity results. In Section~\ref*{sec:4}, several numerical case studies are provided to validate the proposed framework. Finally, discussions and conclusions are presented in Section~\ref*{sec:5}.

\section{Multifidelity hybird solid-porous infill design method} \label{sec:Multifidelity Method}

\begin{figure}[H]
  \centering
  \includegraphics[width=\textwidth]{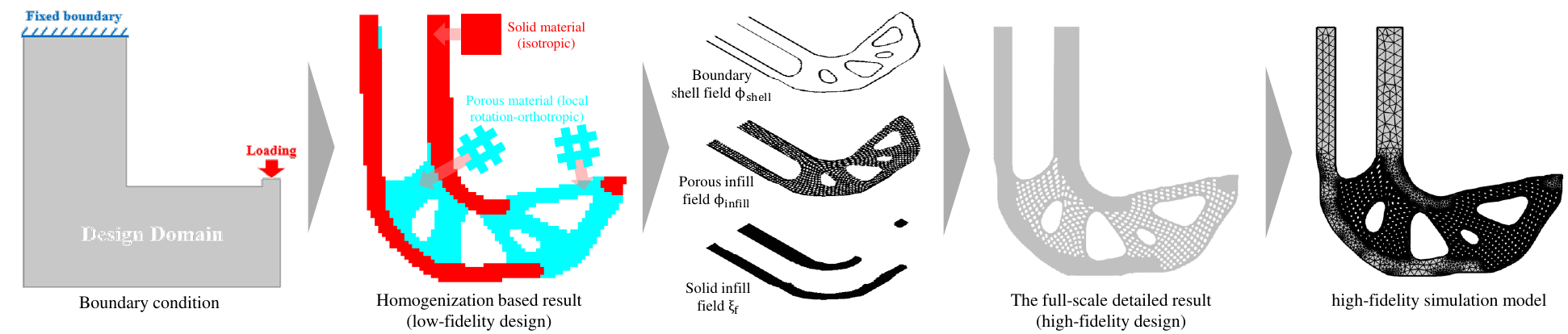}
  \caption{The concept of decomposing the regions contained in a hybrid solid-porous infill structure. }\label{fig:The cpt of hysp struct}
\end{figure}

For a typical hybrid solid-porous infill structure (Figure~\ref{fig:The cpt of hysp struct}), we could decompose it into three regions: shell field, porous infill field, and solid infill field. Shell structure defines the external boundary of the entire structure, porous infill field represents the lattice infill pattern, and solid infill field represents the solid reinforcement regions. By changing the material distribution across these three regions, we can ultimately achieve overall design of the hybrid solid-porous infill structure.

In this section, we introduce the design method for achieving such structures. A key feature of this approach is a multifidelity formulation, which is divided into two main parts: low-fidelity optimization and high-fidelity evaluation. The entire design flowchart is shown in Figure~\ref{fig:The cpt of hysp struct}.

First, for the low-fidelity optimization, we employ a multi-material topology optimization algorithm to preliminarily obtain voxel/pixel-based results. These results consist of three scalar fields, which are considered the low-fidelity outputs. Subsequently, a series of mapping techniques are applied to transform these low-fidelity results into structures with practical geometric significance. Adaptive meshing is then performed on these structures, followed by more accurate Finite Element Analysis (FEA) to evaluate their structural performance.

\subsection{Low fidelity optimization: stress based solid porous infill topology optimization} \label{subsec:2.1}

\begin{figure}[H]
  \centering
  \includegraphics[width=\textwidth]{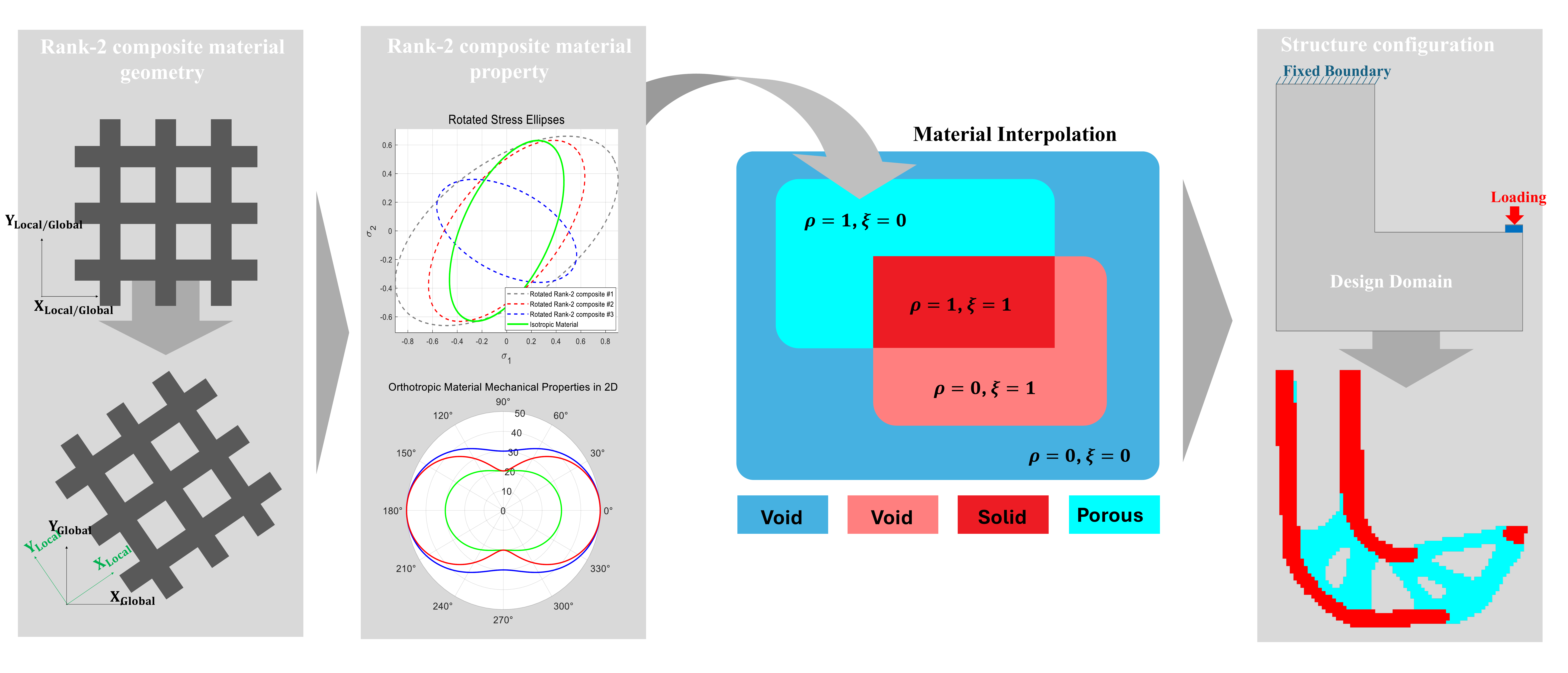}
  \caption{The basic idea of low-fidelity optimization method.}
  \label{fig:LFoptIdea}
\end{figure}

To parameterize solid porous infill composite structures, we adopt the two-field multimaterial design parameterization scheme    \cite{sigmund2001design} that simultaneously characterizes material spatial occupancy and material phase distribution. In this study, a material phase is defined as orthotropic material with a specific orientation (RANK-2 composite material for porous domain shown in Figure~\ref{fig:LFoptIdea}) or isotropic material (solid domain). To better match with the density-based topology optimization framework, two discretized control variable fields (\(\bm{\rho}\) and \(\bm{\xi}\)) are introduced. Within the optimization process, to ensure that the phase fields generated by \(\bm{\rho}\) and \(\bm{\xi}\) are clear and free from single-node connections, we need to seperately perform sequential filtering and projection operations on them:

\begin{equation}
  \begin{aligned}
    \tilde{*} = M(*, R_\text{s}), \quad * = \rho_e, \xi_e
  \end{aligned}
  \label{eq:LFfilter}
\end{equation}
  
The function \(M(*)\) is the smoothing function, and \(R_\text{s}\) is the smoothing radius. It can be any smoothing filter commonly used in topology optimization (such as the PDE filter or convolution filter), for which detailed explanations can be found in \cite{bourdin2001filters}\cite{lazarov2011filters}. Then, the approximated Heaviside projection \cite{wang2011projection} is applied to obtain the 0-1 fields:

\begin{equation} 
  \begin{aligned}
    \bar{*} = \frac{\tanh(\beta_\text{s} \delta_\text{s}) + \tanh(\beta_\text{s} (\tilde{*} - \delta_\text{s}))}{\tanh(\beta_\text{s} \delta_\text{s}) + \tanh(\beta_\text{s} (1 - \delta_\text{s}))}, \quad \tilde{*} = \tilde{\rho}_e, \tilde{\xi}_e
  \end{aligned}
  \label{eq:LFproject}
\end{equation}

where $\beta_\text{s}$ and $\delta_\text{s}$ are the sharpness and the threshold factor. 

\subsubsection{Stiffness interpolation strategy} \label{subsec:2.1.1}

Based on the above two projected control variables, the material stiffness interpolation for mechanical behavior characterization of solid-porous infill composites thus could be expressed by an extended SIMP formulation as follows:

\begin{equation} 
  \mathbb{C}_e(\bar{\rho}_e, \bar{\xi}_e) = 
  \left[ \varepsilon + (1 - \varepsilon)\bar{\rho}_e^{p_\rho} \right]
  \left[ \bar{\xi}_e^{p_\xi} \mathbb{C}^{(\text{iso})} + (1 - \bar{\xi}_e)^{p_\xi} \mathbb{C}_e^{(\text{aniso})} \right]
  \label{eq:extendedSIMP}
\end{equation} 

where \(\mathbb{C}^{(\text{iso})}\) and \(\mathbb{C}_e^{(\text{aniso})}\) are elastic tensors of solid and porous materials, respectively, \(p_\rho\) and \(p_\xi\) are SIMP penalization parameters, and \(\varepsilon\) is a small number to avoid numerical singularity. This study uses \(p_\rho = p_\xi = 3\). It can be seen from the Figure~\ref{fig:LFoptIdea}, when \(\bar{\rho}_e = 1\) and \(\bar{\xi}_e = 1\), \(\mathbb{C}_e\) reduces to \(\mathbb{C}^{(\text{iso})}\); and when \(\bar{\rho}_e = 1\) and \(\bar{\xi}_e = 0\), \(\mathbb{C}_e\) reduces to \(\mathbb{C}_e^{(\text{aniso})}\). Notice that the interpolation in Eq.(\ref{eq:extendedSIMP}) is a popular numerical scheme to facilitate the design optimization rather than computing the homogenized properties, and hence, \(\mathbb{C}_e\) is physically well-defined only when \(\bar{\rho}_e, \bar{\xi}_e \in \{0,1\}\), which is enforced by the projection formulation in Eq.(\ref{eq:LFproject}). In this work, the interpolation in Eq.(\ref{eq:extendedSIMP}) includes different local orientations of the RANK-2 composite material that require stiffness tensor transformations from respective local coordinates to global coordinates. The transformed orthotropic stiffness tensor \(\mathbb{C}_e^{(\text{aniso})}\) in global coordinates is obtained as

\begin{equation} 
  \mathbb{C}_e^{(\text{aniso})} = \mathbb{T}^{-1}(\theta_{e}) \mathbb{C}_0^{(\text{aniso})} \mathbb{T}^{\mathsf{-T}}(\theta_{e})
  \label{eq:GlobalRotate}
\end{equation} 

For 2D plane stress problems used in this study, the matrix forms of the stiffness tensors \(\mathbb{C}_0^{(\text{aniso})}\) and \(\mathbb{C}^{(\text{iso})}\) are given by

\begin{equation} 
  \mathbb{C}_0^{(\text{aniso})} = 
\begin{bmatrix}
\frac{E_{11}}{1-{\nu}_{12}\nu_{21}} & \frac{\nu_{12}E_{22}}{1-\nu_{12}\nu_{21}} & 0 \\
\frac{\nu_{12}E_{22}}{1-\nu_{12}\nu_{21}} & \frac{E_{22}}{1-\nu_{12}\nu_{21}} & 0 \\
0 & 0 & G_{12}
\end{bmatrix},
\quad
\mathbb{C}^{(\text{iso})} = \frac{E_{\text{iso}}}{1 - \nu_{\text{iso}}^2}
\begin{bmatrix}
1 & \nu_{\text{iso}} & 0 \\
\nu_{\text{iso}} & 1 & 0 \\
0 & 0 & \frac{1-\nu_{\text{iso}}}{2}
\end{bmatrix},
\label{eq:ElasticTensor}
\end{equation} 

and the matrix form of the transformation tensor \(\mathbb{T}(\theta_e)\) is given by

\begin{equation} 
\mathbb{T}(\theta_e) = 
\begin{bmatrix}
\cos^2(\theta_e) & \sin^2(\theta_e) & 2\sin(\theta_e)\cos(\theta_e) \\
\sin^2(\theta_e) & \cos^2(\theta_e) & -2\sin(\theta_e)\cos(\theta_e) \\
-\sin(\theta_e)\cos(\theta_e) & \sin(\theta_e)\cos(\theta_e) & \cos^2(\theta_e) - \sin^2(\theta_e)
\end{bmatrix},
\label{eq:RotateMatrix}
\end{equation} 

where \(E_{11}, E_{22}, \nu_{12}, \text{ and } G_{12}\) are elastic moduli along the local direction, elastic modulus perpendicular to the local direction, Poisson's ratio with respect to the local direction, and shear modulus, respectively, for the orthotropic composite material. \(E_{\text{iso}}\) and \(\nu_{\text{iso}}\) are elastic modulus and Poisson's ratio, respectively, for the isotropic material. Generally, the local direction \(\bm{\theta}\) is assumed along with the principle stress direction. For a 2D plane stress problem, the elemental second-order stress tensor $\bm{\mathrm{\sigma}}$ is:
\begin{equation}
  \bm{\mathrm{\sigma}}_e = \begin{bmatrix}
    {\sigma}_{xx,e} & {\sigma}_{xy,e} \\
    {\sigma}_{yx,e} & {\sigma}_{yy,e}
  \end{bmatrix} 
  \label{eq:StressTensor}
\end{equation}

where $\sigma_{xx}, \sigma_{yy}$ are the normal stress along $x$ and $y$ direction respectively (the direction shown in Figure~\ref{fig:LFoptIdea}), and $\sigma_{xy}$ is the shear stress within the $x$-$y$ plane. For the principal direction \(\bm{\theta}\), we could directly solve the following equation to obtain:

\begin{equation}
\tan(\theta_e) = \frac{\sigma_{xy,e}}{\sigma_{xx,e} - \sigma_{yy,e}}
\label{eq:PSLdirection}
\end{equation}

Generally, the $\theta_e$ calculated from Eq.(\ref{eq:PSLdirection}) ranges from $-\pi$ to $\pi$. To better suit the subsequent optimization, $\theta_e$ is further transformed into a 0-1 scalar form as follows:

\begin{equation} 
  \bar{\theta_e} = \frac{\theta_e + \pi}{2\pi}
  \label{eq:PSLdirection2}
\end{equation}

\subsubsection{Stress state interpolation strategy} \label{subsec:2.1.2}
To better achieve multi-material stress minimization optimization, in this context, we introduce a term $\bm{\sigma}_{\text{stress}}$ called ``elemental stress state," which does not refer to the actual element stress values but rather an equivalent stress index value. This index value is computed by normalizing with respect to the yield stress, based on different stress criteria. For solid isotropic materials, we employ the von Mises stress criterion under the plane stress assumption, which can be expressed using the following normalized equation:

\begin{equation} 
  \sigma_{\text{von-Mises},e} = f_{\text{vm}}(\bm{\sigma}_e) = \left( \frac{\sigma_{xx,e}}{\sigma_\text{S}} \right)^2 + \left( \frac{\sigma_{yy,e}}{\sigma_\text{S}} \right)^2 - \frac{\sigma_{xx,e} \sigma_{yy,e}}{\sigma_\text{S}^2} + \left( \frac{3\sigma_{xy,e}}{\sigma_\text{S}} \right)^2
  \label{eq:VMstressState}
\end{equation}

where \(\sigma_\text{S}\) is the yield strength of the isotropic material. 

For the orthotropic materials in the porous regions, considering the characteristics of the microstructure, we assume identical tensile and compressive properties. Therefore, the plane stress based Tsai-Hill stress criterion is used:

\begin{equation}
  \sigma_{\text{Tsai-Hill,e}} = f_{\text{th}}(\bm{\sigma}_e,\theta_e) = \left( \frac{\sigma_{11,e}}{\sigma_{\text{X}}} \right)^2 + \left( \frac{\sigma_{22,e}}{\sigma_{\text{Y}}} \right)^2 - \frac{\sigma_{11,e} \sigma_{22,e}}{\sigma_{\text{X}} \sigma_{\text{Y}}} + \left( \frac{\sigma_{12,e}}{\sigma_{\text{XY}}} \right)^2
  \label{eq:THstressState}
\end{equation} 

where \(\sigma_{\text{X}}\), \(\sigma_{\text{Y}}\), and \(\sigma_{\text{XY}}\) represent the yield strengths of the anisotropic material in different directions. \({\sigma}_{11}\), \({\sigma}_{22}\), and \({\sigma}_{12}\) are the in-plane stress components in the local coordinate system, obtained by transforming the global coordinate in-plane stress components as:

\begin{equation} 
  \begin{bmatrix}
    \sigma_{11,e} \\
    \sigma_{22,e} \\
    \sigma_{12,e}
    \end{bmatrix}
    =
    \mathbb{T}(\theta_e)
    \begin{bmatrix}
    \sigma_{xx,e} \\
    \sigma_{yy,e} \\
    \sigma_{xy,e}
    \end{bmatrix}
  \label{eq:Local2global}
\end{equation} 

where \(\mathbb{T}(*)\) is the transformation tensor in Eq.(\ref{eq:RotateMatrix}). Using the interpolation in Eq.(\ref{eq:VMstressState}) and Eq.(\ref{eq:THstressState}), we propose the interpolated yield function for composite structures with orthotropic Tsai-Hill and isotropic von Mises candidate materials as:

\begin{equation} 
  \sigma_{\text{stress},e}(\bar{\xi}_e, \bm{\mathrm{\sigma}}_e,\theta_e) = \bar{\xi}_e^{p_\xi} f_\text{vm}(\bm{\mathrm{\sigma}}_e) + (1 - \bar{\xi}_e)^{p_\xi} f_\text{th}(\bm{\mathrm{\sigma}}_e,\theta_e)
  \label{eq:13}
\end{equation}  

With the introduction of parameters such as \(\sigma_\text{S}\), \(\sigma_\text{X}\), \(\sigma_\text{Y}\), and \(\sigma_\text{XY}\), this formulation can accurately capture the stress responses of the two different materials. Clearly, for \(\sigma_{\text{stress},e}\), a smaller value indicates a better overall stress state for the structure.

\subsubsection{Material volume interpolation strategy} \label{subsec:2.1.3}
Recalling the expression in Eq.(\ref{eq:extendedSIMP}), the \(\bar{\rho}\) is used to determine the base material distributions:
\begin{equation} 
  v_e^\text{base} = \bar{\rho}_e
  \label{eq:14}
\end{equation} 
and \(\bar{\xi}\) is used to determine the solid material distribution. For the elemental volume of solid material, we could express as:
\begin{equation} 
  v_e^\text{solid} = \bar{\xi}_e\bar{\rho}_e
  \label{eq:15}
\end{equation} 

Additionally, we seperately constrain the volume fraction for the solid part indicated by the field in Eq.(\ref{eq:15}) and the overall base structure \(\bar{\bm{\rho}}\) to indirectly control the usage ratio of solid material and porous infill material.

\subsubsection{Low fidelity optimization formulation} \label{subsec:2.1.4}

Since the primary purpose of solid infill is to enhance the local strength of the structure to counteract local stress concentrations, we take the maximum stress state as the optimization objective. To ensure that this objective function is stably differentiable, we use a P-norm aggregation function to approximate the maximum value:

\begin{equation} 
\max\left(\forall \sigma_{\text{stress},e}\right) \approx \sigma_{\text{PN}} = \left( \sum_{e=1}^{N_\text{c}} \left( \sigma_{\text{stress},e} \right)^P \right)^{\frac{1}{P}}
\label{eq:16}
\end{equation}

where \(\sigma_{\text{PN}}\) is the global P-norm measure, \(P\) is the aggregation parameter, and \(N_\text{c}\) denotes the number of elements for \(\bm{\rho}\) and \(\bm{\xi}\). Note that the P-norm approaches the maximum stress as \(P \to \infty\). A large \(P\) value tends to make the stress-constrained problem ill-conditioned. Relatively small \(P\) values are preferred in practice to ensure convergence stability, which, however, introduces a gap between the exact and approximated maximum stresses. To better approximate the maximum stress without excessively increasing the \(P\) value, the global P-norm stress measure is iteratively corrected through:

\begin{equation} 
\bar{\sigma}_{\text{PN}} = c^{(L)} \sigma_{\text{PN}}
\label{eq:17}
\end{equation}

where \(c^{(L)}\) is the correction parameter at the \(L^\text{th}\) iteration \((L > 1)\) that reflects the ratio of the maximum von Mises stress to the P-norm stress from the current iteration. Note that the change in \(c^{(L)}\) can be abrupt if the correction is based solely on the history-independent stress ratio, leading to oscillations and instability in convergence. To mitigate this issue, a parameter \(\alpha \, (\alpha \in (0,1])\) is introduced to restrict the variation between \(c^{L}\) and \(c^{(L-1)}\), as shown below \cite{yang2018stress}:

\begin{equation} 
  c^{(L-1)} = \alpha \frac{\max \left( \sigma_{\text{stress},e} \right)}{\sigma_{\text{PN}}} + (1 - \alpha) c^{(L-1)}
\label{eq:18}
\end{equation}

In this work, \(\alpha = 0.5\) is adopted for all iterations and \(c^{(0)} = 1\) is used.

The topology optimization problem is solved in nested form, by successive minimizations with respect to design variables \(\bm{\rho}\), \(\bm{\xi}\), and \(\bm{\theta}\); where for each design iteration, the equilibrium equations are satisfied by FEA. As is shown by Pedersen \cite{pedersen1989optimal}, the optimal orientation of an orthotropic composite coincides with the principal stress directions, hence \(\bm{\theta}\) is aligned accordingly for each minimization step and not directly optimized. Subsequently, design vectors \(\bm{\rho}\) and \(\bm{\xi}\) are updated at each minimization step based on their gradients using MMA \cite{svanberg1987method}. The discretized optimization problem can thus be written as:

\begin{equation}
    \begin{array}{ll}
        \text{find:} & \bm{\rho}, \bm{\xi} \\
        \text{minimize:} & \sigma_\text{PN}\\
        \text{subject to:} & \begin{cases}
          \mathbf{K}\mathbf{U} = \mathbf{F} & \text{(19-1)} \\
        \frac{\sum_{e=1}^{N_\text{c}} v_e^\text{base}}{N_\text{c}} \leq V_\text{b} & \text{(19-2)} \\
        \frac{\sum_{e=1}^{N_\text{c}} v_e^\text{solid}}{{N_\text{c}} V_\text{b}} \leq V_\text{s} & \text{(19-3)} \\
        0 \leq \rho_e,\xi_e \leq 1, \quad \forall e & \text{(19-4)}
        \end{cases}
    \end{array} \label{eq:LFmodel}
\end{equation}

In this context, Eq.(19-1) corresponds to the control equation system of the physical problem. Eq.(19-2) and Eq.(19-3) represents the first and second volume fraction constraint of the structure, where \(V_\text{b}\) and \(V_\text{s}\) are the user-defined upper limit for the mass ratio for base and solid material respectively.

\begin{figure}[H]
  \centering
  \begin{subfigure}[b]{1\textwidth}
    \centering
    \includegraphics[width=\textwidth]{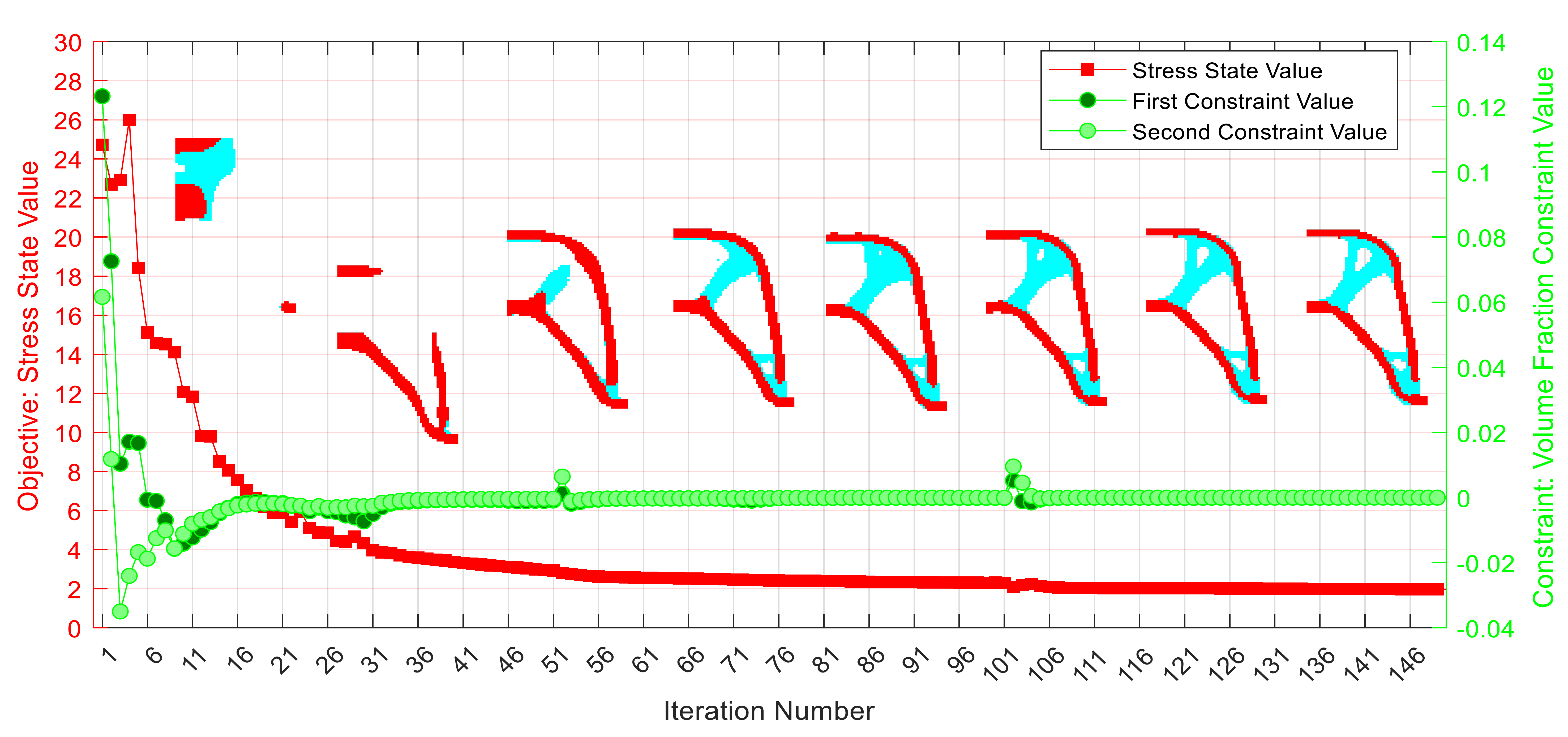} 
    \caption{}
  \end{subfigure}

  \vspace{0.5cm} 
  \begin{subfigure}[b]{1\textwidth}
    \centering
    \includegraphics[width=\textwidth]{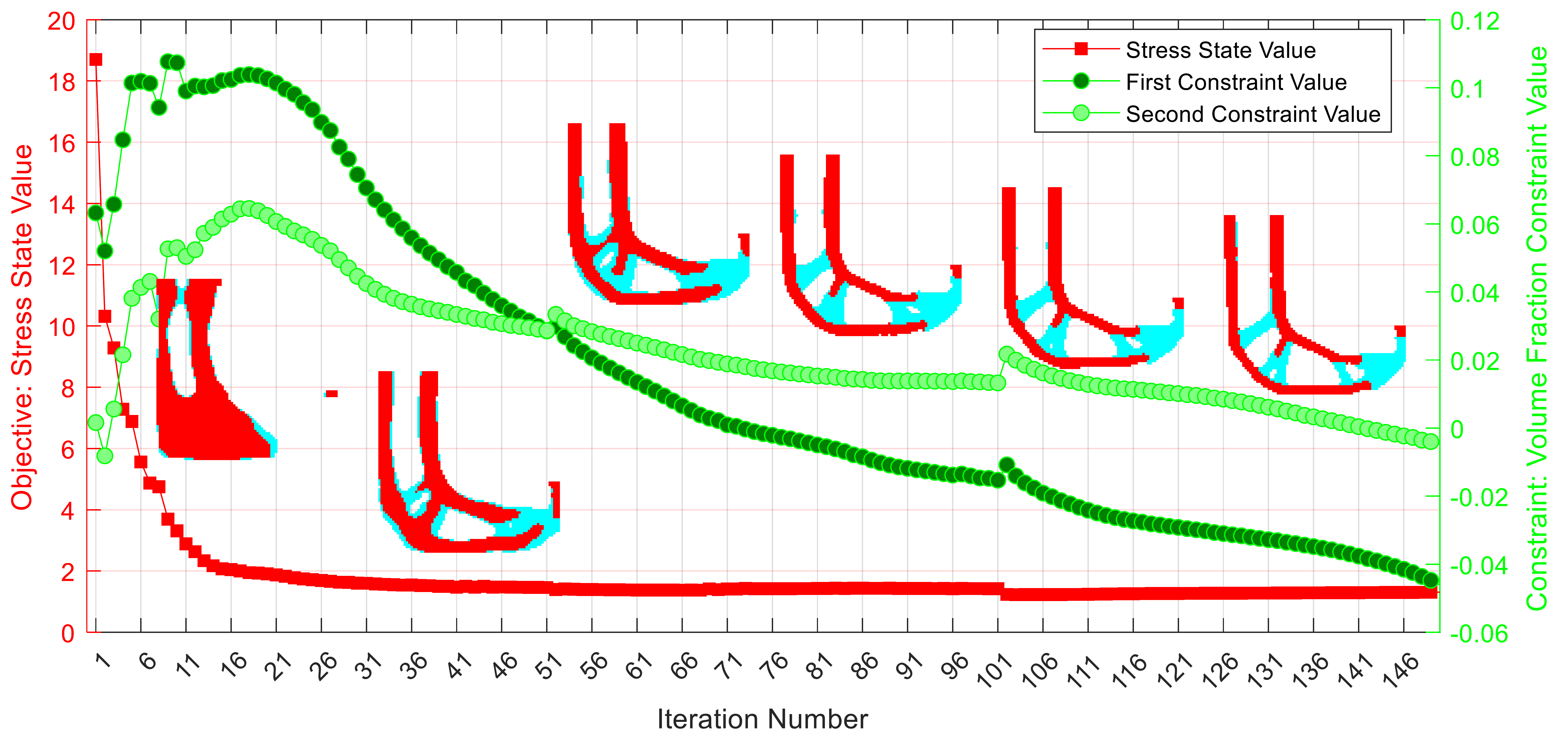} 
    \caption{}
  \end{subfigure}

  \caption{The optimization process for low-fidelity optimization: (a) symmetric tension beam case, and (b) L-bracket beam case.}
  \label{fig:The LFopt process}
\end{figure}

Apparently, Eq.(\ref{eq:LFmodel}) represents a typical multiple-material, maximum stress minimization-based topology optimization problem, which is well-studied and can be solved efficiently. Specific details regarding the optimization, such as adjoint sensitivity analysis, are omitted here. For more information, refer to the following classic papers on stress-based topology optimization \cite{le2010stress}\cite{holmberg2013stress}\cite{yang2018stress} or multi-material based stress topology optimization \cite{xu2021stress}\cite{kundu2022multimaterial}\cite{kundu2023stress}\cite{nguyen2023improving}\cite{nguyen2024design}. Figure~\ref{fig:The LFopt process} illustrates the optimization process for two cases of the proposed solid infill design. The effectiveness of this low-fidelity optimization will be discussed in subsequent section.

\subsection{High-fidelity evaluation: full-scale structure realization through de-homogenization} \label{subsec:2.2}
The process of converting homogenization-based results into macroscopic results with practical geometric significance is referred to as de-homogenization. In this paper, the implementation of this process involves three steps. First, the abstract density elements representing porous regions in the low-fidelity optimization results are transformed into macroscopic truss-like structures using a mapping method based on the wave projection function. Second, the solid-filled regions are correspondingly transformed to form shell structures, which are then combined with the macroscopic porous structure generated in the first step. This results in a realistic hybrid solid-porous infill structure with an actual geometric shape, exhibiting performance similar to the homogenization-based results. Finally, high-fidelity simulation analysis is conducted on the resulting structure. The following two subsections will introduce these processes in detail.

\subsubsection{Wave projection function method} \label{subsec:2.1.1}
\begin{figure}[H]
  \centering
  \includegraphics[width=\textwidth]{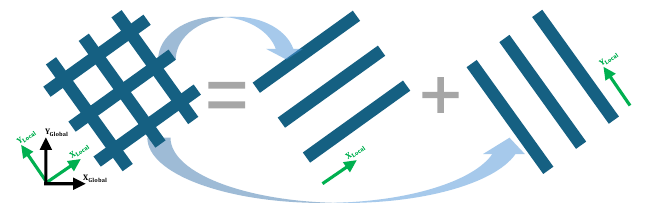}
  \caption{Geometry pattern for rank-2 composite material.}\label{fig:RANK-2 composite}
\end{figure}

Almost any type of periodic microstructure can be represented by a complex exponential fourier series with spatially varying parameters. This approach allows for the projection of a complex microstructure on a fine scale while maintaining a smooth and continuous lattice. For the RANK-2 composite material adopted in this work, which consists of two stripe-like spatial features, to achieve the infill. As shown in Figure~\ref{fig:RANK-2 composite}, it is simple enough to be represented by just two independent fields. The first scalar field, \(\Gamma^{(1)}\), describes the part of the unit cell aligned with the local \(x\)-direction, while the second scalar field, \(\Gamma^{(2)}\), describes the part aligned with the local \(y\)-direction. Simultaneously, we aim for the lattice infill to be reasonably distributed at specific angles in space, thereby achieving weight reduction while also demonstrating certain mechanical advantages. These two fields are derived independently from two gradient vector fields, \(\bm{K}^{(1)}\) and \(\bm{K}^{(2)}\), to achieve the desired distribution.

The wave projection method is used to build the relationship between the vector field $\bm{K}$ and the scalar field $\Gamma$. Let us start by assuming a scalar field $\Phi(\bm{r})$ expressed by the following equation:

\begin{equation} 
  \Phi(\bm{r}) = \bm{K}(\bm{r}) \cdot \bm{r} \label{eq:20}
\end{equation} 

where $\bm{r}$ is the position vector. Then, substituting the phase field $\Phi(\bm{r})$ into the wave function, we could arrive:

\begin{equation} 
  \Gamma(\bm{r}) = \frac{1}{2} + \frac{1}{2} \cos(2\pi \Phi(\bm{r})) \label{eq:21}
\end{equation} 

\begin{figure}[H]
  \centering
  \includegraphics[width=\textwidth]{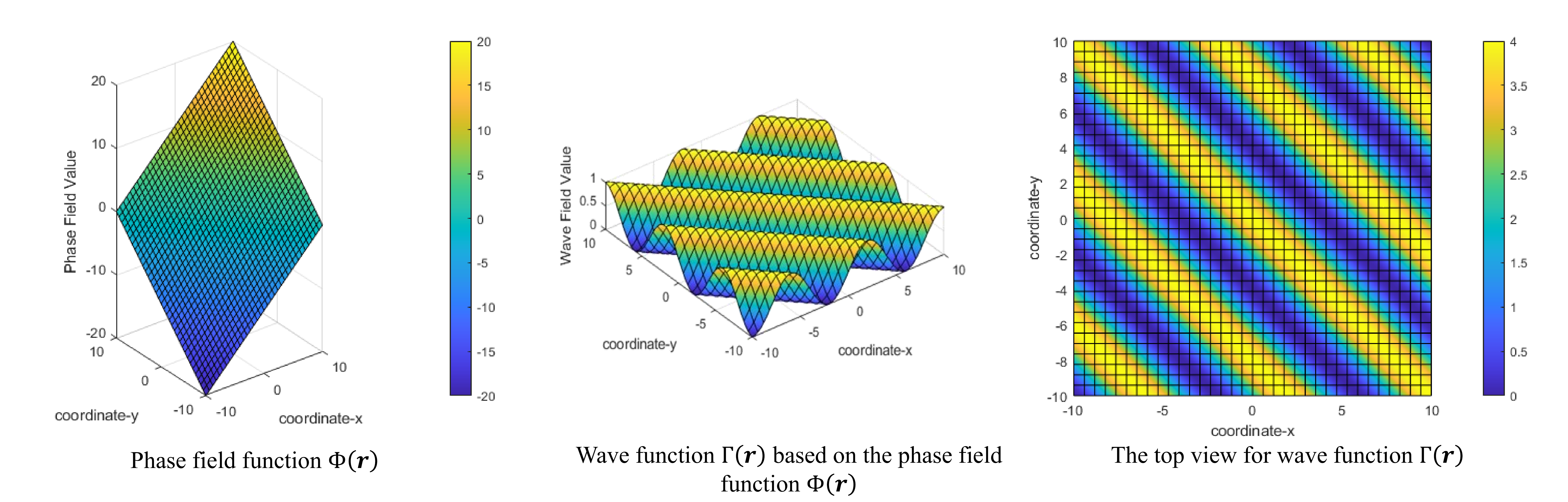}
  \caption{The schematic diagram for the phase field $\Phi(\bm{r})$, the wave function field $\Gamma(\bm{r})$, and their relationship}
  \label{fig:the phase field}
\end{figure}

As shown in Figure~\ref{fig:the phase field}, Eq.(\ref{eq:20}) could make a coordinate-based continuous scalar field $\Phi(\bm{r})$ to be transformed into a spatially periodic scalar field $\Gamma(\bm{r})$, and the peak region of the field $\Gamma(\bm{r})$ is consistent with the distribution of the isoline of the phase field $\Phi(\bm{r})$. Therefore, $\Gamma(\bm{r})$ can project the vector field $\bm{K}$ to result in the formation of stripe-like spatial features.

To obtain the infill pattern mapping with a specific distribution within a design domain, we substitute Eq.(\ref{eq:20}) into Eq.(\ref{eq:21}) to form the wave projection function \(\Gamma(\bm{r})\) as follows:

\begin{equation} 
  \Gamma(\bm{r}) = \frac{1}{2} + \frac{1}{2} \cos(P_\text{d} \bm{K}(\bm{r}) \cdot \bm{r}) \label{eq:22}
\end{equation} 

where $P_\text{d}=\frac{\pi}{d}$, and $d$ is used to control the periodic of the wave function, and

\begin{equation} 
  \bm{K} = 
  \left[
    -\sin\left(\theta(\bm{r})\right),
    \cos\left(\theta(\bm{r})\right)
    \right]^\mathsf{T}
  \label{eq:23}
\end{equation} 

where \(\theta\) is the direction-based field that indicates \(\bm{K}\). For a given directional field \(\theta\), a corresponding cosine function field \(\Gamma\) can be generated, leading to the formation of the corresponding spatial structural features. However, the field \(\Gamma\) cannot be obtained directly through Eq.(\ref{eq:22}) when the propagation direction \(\bm{K}\) is spatially variant. Therefore, the phase field \(\Phi(\bm{r})\) needs to be specifically established to match the spatial variation of \(\bm{K}\). Based on Eq.(\ref{eq:20}), clearly, the vector field $\bm{K}$ is the gradient of the phase field $\Phi(\bm{r})$, and it is perpendicular to the isoline of the phase field $\Phi(\bm{r})$, which can be described more generall by the following equation:

\begin{equation} 
  \nabla \Phi(\bm{r}) = \bm{K}(\bm{r})  \label{eq:24}
\end{equation} 

To solve the above equation, various numerical approaches can be employed. In this work, a least-squares method is applied. The formulation of the least-squares approach is expressed as follows:

\begin{equation}
  \begin{split} 
    \min_{\Phi(\bm{r})} L(\Phi(\bm{r})) &= \int_{\Omega} \|\nabla \Phi(\bm{r}) - \bm{K}(\bm{r})\|^2 \, d\Omega \\
    \quad \text{s.t.} &\quad \nabla \Phi(\bm{r}) \cdot \bm{K}_*(\bm{r}) = 0 \label{eq:25}
  \end{split}
\end{equation} 

where
\begin{equation} 
  \bm{K}_*(\bm{r}) = \left[
    -\sin\left(\theta(\bm{r})+\frac{\pi}{2}\right),
    \cos\left( \theta(\bm{r})+\frac{\pi}{2}\right)
    \right]^\mathsf{T}
    \label{eq:26}
\end{equation} 

\(\bm{K}_*(\bm{r})\) is the result of rotating the vector field of the given angle field \(\theta\) by 90 degrees. The reason for adding an additional constraint in Eq.(\ref{eq:25}) is to prevent the significant angular discrepancy between \(\Phi(\bm{r})\) and \(\bm{K}(\bm{r})\) in regions where the angle changes abruptly.

\subsubsection{Explicit geometry feature control} \label{subsec:2.2.2}
From the previous discussion, we know that the peak regions of the cosine function \(\Gamma(\bm{r})\) align with the distribution of the isolines of the phase function \(\Phi(\bm{r})\) obtained from a given directional field \(\theta\). To scale the microstructure to the size of bilinear elements, we assume \(\Gamma\) has a discretized form \(\bm{\Gamma}\) in space, and a series of graphical operations needs to be performed for \(\bm{\Gamma}\). 
Firstly, the \(\bm{\Gamma}\) is converted into a 0-1 field \(\bar{\bm{\Gamma}}\):

\begin{equation} 
  \bar{\Gamma}_e =
\begin{cases} 
      0, & \text{if } \Gamma_e < 0.5 \\
      1, & \text{if } \Gamma_e \geq 0.5
\end{cases}
\label{eq:27}
\end{equation}

To obtain a fin structure with uniform width control, we will perform skeleton extraction on $\bar{\bm{\Gamma}}$ to create a skeleton density field, defined as $\bar{\bm{\Gamma}}_{\text{skel}}$.

\begin{equation} 
    \bar{\bm{\Gamma}}_{\text{skel}} = \text{bwskel}(\bar{\bm{\Gamma}}) \label{eq:28}
\end{equation} 

Here, \(\text{bwskel}(*)\) reduces all objects in \(\bar{\bm{\Gamma}}\) to 1-pixel wide curved lines without changing the essential structure of the image. This process extracts the centerline while preserving the topology of the objects. To ensure manufacturability and control the width, a further dilation process is applied to \(\bar{\bm{\Gamma}}_{\text{skel}}\):

\begin{equation} 
    \bar{\bm{\Gamma}}_{\text{d}} = \text{imdilate}(\bar{\bm{\Gamma}}_{\text{skel}}, \mathbb{{E}}) \label{eq:29}
\end{equation} 

\begin{figure}[H]
  \centering
    \includegraphics[width=\textwidth]{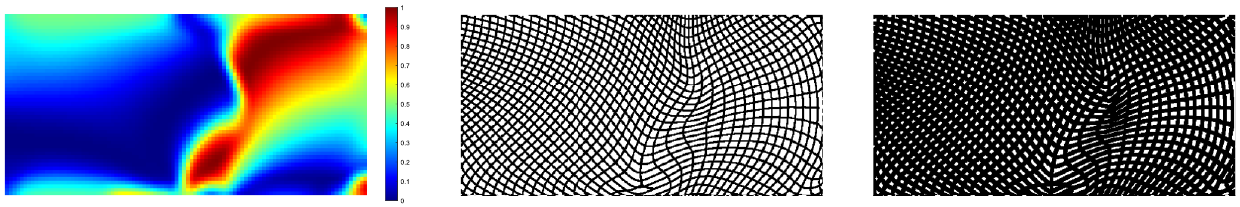} 

  \caption{The given field \(\bm{\theta}\) with random boundary condition (left), the corresponding lattice pattern \(\bar{\bm{\Gamma}}_{\text{d}}\) with \(d=8\) and \(N_{\text{dilate}}=3\) (middle), and the corresponding lattice pattern \(\bar{\bm{\Gamma}}_{\text{d}}\) with \(d=8\) and \(N_{\text{dilate}}=6\) (right).}
  \label{fig:lattice infill patterns}
\end{figure}

The \(\text{imdilate}(*)\) function performs the dilation operation. By adjusting the size (\(N_{\text{dilate}}\)) of the morphological kernel \(\mathbb{{E}}\), the width of the structure within the field \(\bar{\bm{\Gamma}}_{\text{d}}\) can be controlled. Figure~ \ref{fig:lattice infill patterns} shows two lattice infill patterns with different \(\bm{\theta}\) inputs and different geometry parameter controls.

\subsubsection{Final structure realization and simulation} \label{subsec:2.2.3}
Based on the content from Section~\ref*{subsec:2.2.2}, we have a general understanding of how to transform the low-fidelity data into the necessary structure for a hybrid solid-porous infill geometry, enabling a more precise and comprehensive analysis of the problem at hand. However, to make this process more controllable and better integrate with the results from Section~\ref*{subsec:2.1}, additional operations are required. The process is structured into six main steps:

\textbf{Step 1:} First, we need to extract the three relevant design domains from the low-fidelity optimization—namely, the base material field, the solid material field, and the principal stress direction field—as the low-fidelity results. These discretized fields are denoted as \(\bm{\rho}_\text{c}\), \(\bm{\xi}_\text{c}\), and \(\bm{\theta}_\text{c}\), respectively. Here, the subscript \(\text{c}\) indicates that their dimensions correspond to the number of design variables \(N_\text{c}\) in the low-fidelity optimization. Subsequently, these fields are refined into new fields, \(\bm{\rho}_\text{f}\), \(\bm{\xi}_\text{f}\), and \(\bm{\theta}_\text{f}\), each with the refined dimension \(N_\text{f}\), through linear interpolation with the scaling factor \(\lambda\).

\textbf{Step 2:} Subsequently, shell extraction is performed based on the obtained \(\bm{\rho}_\text{f}\), resulting in $\bm{\phi}_{\text{shell}}$. 

\textbf{Step 3:} The operations introduced in Subsection \ref*{subsec:2.2} are applied to \(\bm{\theta}_\text{f}\) to obtain the lattice infill pattern field \(\bar{\bm{\Gamma}}_{\text{d}}\).

\textbf{Step 4:} For the obtained refined fields $\bm{\rho}_\text{f}$, \(\bm{\xi}_\text{f}\), $\bm{\phi}_{\text{shell}}$, and \(\bar{\bm{\Gamma}}_{\text{d}}\), a series of operations are performed to generate the final structure field $\bm{\phi}$:

\begin{equation} 
  \phi_{\text{infill},e} = \rho_{\text{f},e}\bar{\Gamma}_{\text{d},e} + \phi_{\text{shell},e}
  \label{eq:30}
\end{equation}

This operation is used to determine the lattice infill domain \(\bm{\phi}_{\text{infill}}\) within the given structure. Then, a 0-1 projection is applied as follows:

\begin{equation} 
  \phi_{\text{infill},e}' =
\begin{cases} 
      0, & \text{if } \phi_{\text{infill},e} < 0.5 \\
      1, & \text{if } \phi_{\text{infill},e} \geq 0.5
\end{cases}
\label{eq:31}
\end{equation}

The Boolean addition operation is then applied at each element to achieve the solid infill:

\begin{equation} 
  \phi_e = \min(\xi_{\text{f},e}, \phi_{\text{infill},e}')
  \label{eq:32}
\end{equation}

The obtained structure field \(\bm{\phi}\) will be further smoothed with a smoothing radius \(R_\text{f}\):

\begin{equation} 
\tilde{\phi_e} = M(\phi_e, R_\text{f})
\label{eq:33}
\end{equation}

\textbf{Step 5:} the explicit boundary will be extracted to form the smoothed field \(\tilde{\bm{\phi}}\). This operation can be performed using contour extraction function:

\begin{equation} 
  \mathbb{{Y}} = \text{contourc}(\tilde{\bm{\phi}}, \eta_\text{t})
\label{eq:34}
\end{equation}

where \(\mathbb{{Y}}\) is a new data structure containing the ordered coordinates of multiple points that represent polylines showing the explicit boundary. Here, \(\eta_\text{t}\) is the height of the extracted contour line. Since \(0 < \tilde{\phi}_e < 1\), \(\eta_\text{t}\) can take any value between 0 and 1. Clearly, the values of \(R_\text{f}\) and \(\eta_\text{t}\) influence the final stress analysis results. However, because this work optimizes based on a fixed standard value, the sensitivity of the stress results is outside the scope of this article; interested readers can refer to the relevant literature. In this work, we set \(R_\text{f} = 8\) and \(\eta_\text{t} = 0.5\). Finally, the resulting \(\mathbb{{Y}}\) can be saved as a .dxf file, a standard CAD format readable by most CAD softwares.

\textbf{Step 6:} Based on the geometric model (\texttt{.dxf} file) obtained from Step 5, adaptive meshing is performed (this step can also be done easily in commercial CAE software or open-source meshing software). Under the given boundary conditions and physical problem setup, linear elastic FEA is conducted on the mesh in this work, focusing specifically on the maximum von Mises stress (we assume that the final detailed design will be fabricated with isotropic material). This metric is collectively recorded as the objective function values for the further optimization process.

The entire procedure is illustrated in the Figure~\ref{fig:HFprocess}:

\begin{figure}[H]
  \centering
  \includegraphics[width=\textwidth]{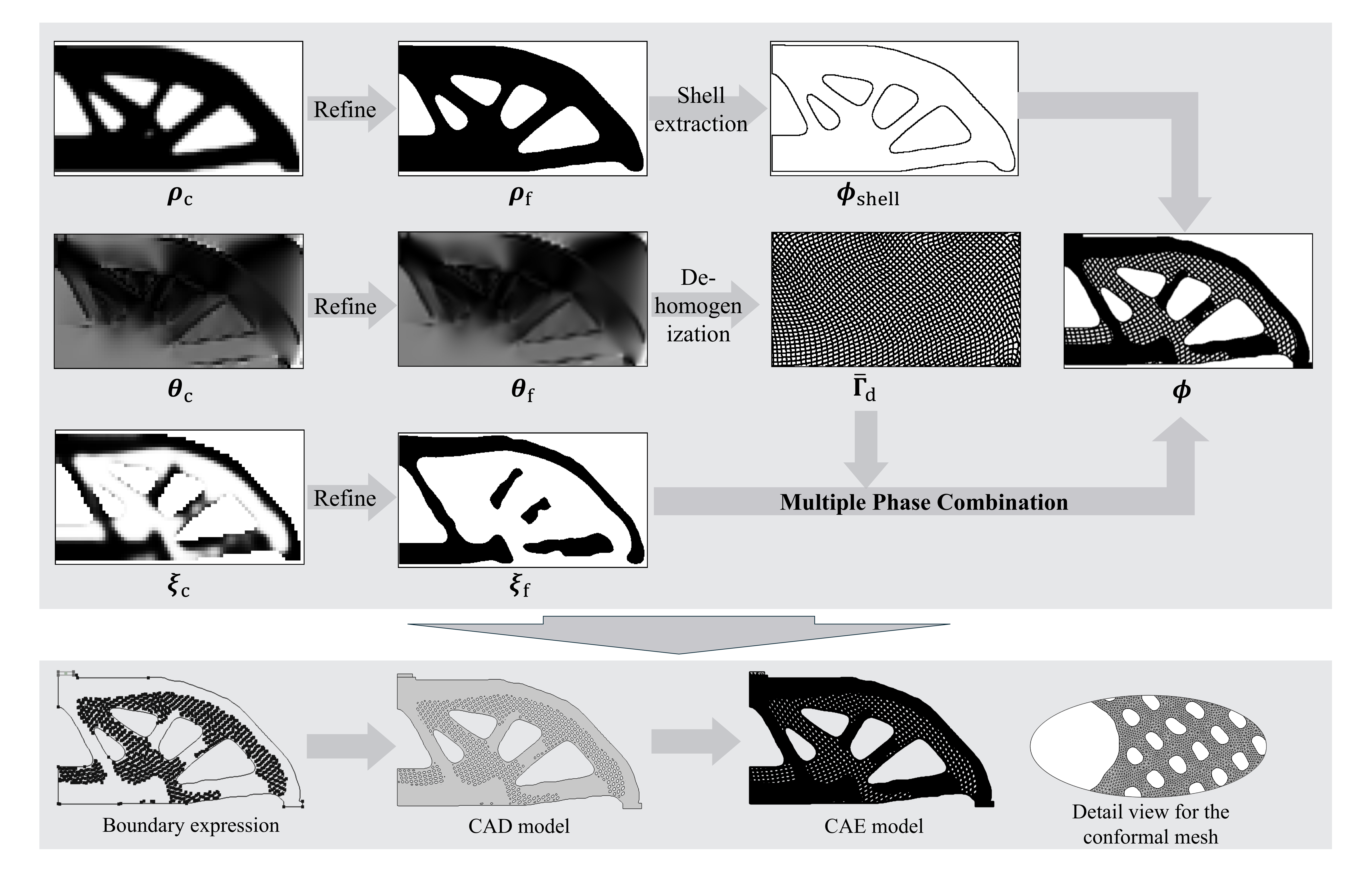}
  \caption{The procedure for obtaining the detail result used for high-fidelity evaluation.} 
  \label{fig:HFprocess}
\end{figure}

\section{Multifidelity optimization framework} \label{sec:3}
\subsection{Optimization model formulation} \label{subsec:3.1}
Recalling the previous content, converting low-fidelity data into high-fidelity results involves numerous complex operations, making it challenging to accurately and efficiently determine the sensitivities of these processes. Moreover, high-fidelity results contain significantly more detailed information than low-fidelity data, and this disparity in information capacity results in highly nonlinear optimization. To address this issue, we use data-driven MFTD method, to execute the rest optimization process.

Similar to traditional EA, population diversity significantly impacts the final optimization results in data-driven MFTD. In the context of structural geometry optimization problems, ensuring population diversity is crucial. The most effective approach to achieve this is by varying the mass ratio of the final optimized structures, which leads to diverse geometric shapes. However, in structural optimization problems focused on strength, it is evident that a higher mass ratio generally results in improved structural performance. If structural performance alone is prioritized, the algorithm may prematurely converge toward higher mass ratio designs, reducing diversity within the population. To address this issue, our algorithm employs a multi-objective optimization approach, considering both the structural mass ratio and performance. By balancing the trade-off between mass ratio and structural performance, we maintain population diversity throughout each optimization iteration, ultimately obtaining a series of diverse design results.

In conclusion, our optimization model can be represented as follows:

\begin{equation}
  \begin{array}{ll}
      \text{find:} & \mathcal{P} = \left\{ \mathbf{x}_1, \mathbf{x}_2, \ldots, \mathbf{x}_\text{m} \right\} \\
      \text{minimize:} & \mathbf{O}(\mathbf{x}_i) = \left[G_\text{vf}(\mathbf{x}_i), G_\text{opt}(\mathbf{x}_i)\right], \quad i = 1, 2, \ldots, \text{m} \\
      \text{subject to:} & \mathbf{x}_i \in \Omega_\text{c}\\
  \end{array} 
  \label{eq:35}
\end{equation}

\(\mathbf{x}_i = \left\{ \mathbf{x}_{1}, \mathbf{x}_{2}, \ldots, \mathbf{x}_\text{m} \right\}\) represents the decision variable vector (the low-fidelity data) of the  individual, where \(\text{m}\) is the number of decision variables. \(G_\text{vf}\) is the mass ratio for the high-fidelity results, \(G_\text{opt}\) is the maximum von Mises stress for the high-fidelity results, and \(\Omega_\text{c}\) denotes the valid design domain, which will be explained in subsequent sections.

\subsection{Data-driven multifidelity topology design framework} \label{subsec:3.2}

\begin{figure}[H]
  \centering
  \includegraphics[width=\textwidth]{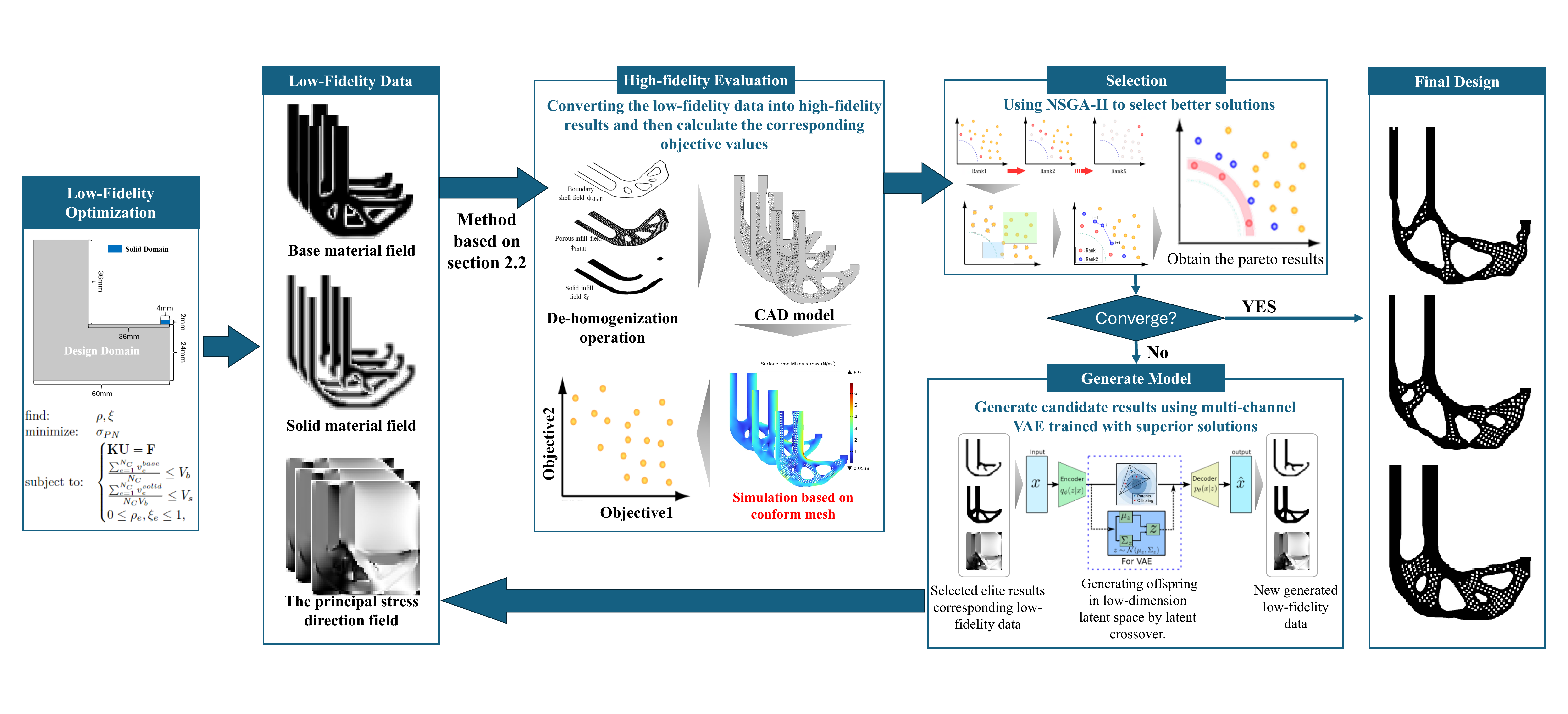}
  \caption{The schematic illustration of data-driven multifidelity topology design.}  
  \label{fig:MFTDflowchat}
\end{figure}

The procedures of data-driven MFTD are shown in Figure~\ref{fig:MFTDflowchat}. Each step is briefly described below. Note that this paper omits a lot of formula details proposed in the original framework for simplicity and the detailed content can be found in the original paper \cite{yaji2022data}\cite{yamasaki2021data}. 

\subsubsection{Low-fidelity Optimization} \label{subsec:3.2.1}

Within the proposed data-driven MFTD framework, the low-fidelity optimization introduced in Section~\ref*{subsec:2.1} is performed only once before the first optimization iteration. During this stage, it is executed multiple times based on the desired number of decision variables, \(\text{m}\), to generate the initial samples for the subsequent optimization algorithm.

It is important to note that many simplifications are applied during the low-fidelity optimization process. As a result, the solutions generated from low-fidelity optimization may not always be highly accurate. Nevertheless, its primary function is to provide a sufficient number of ``relatively reliable" initial solutions, which are continuously refined and updated in subsequent optimization iterations to yield more accurate and reasonable results. This aspect will be further elaborated in a later section.

\subsubsection{High-fidelity Evaluation} \label{subsec:3.2.2}

The operation introduced in Subsection \ref*{subsec:2.2} is conducted for each effective low-fidelity data point to obtain the corresponding high-fidelity results. Subsequently, the detailed mass ratio \(G_\text{vf}\) and structural performance index \(G_\text{opt}\) are calculated.

It is worth noting that anomalies may occasionally occur during this process, such as when the low-fidelity data quality is too poor to be processed by high-fidelity simulation, or when the objective function values are excessively high. In such cases, we exclude these anomalous candidates directly to avoid their impact on the algorithm.

\subsubsection{Selection} \label{subsec:3.2.3}

Based on the high-fidelity evaluation results, superior candidates, known as elite solutions, are selected from the candidate pool using the elite strategy. In this study, the NSGA-II, a multi-objective algorithm, is employed to rank and select candidates based on the Pareto dominance relation in the objective space. The current best solutions are named as Rank-1 results, i.e., no other solutions dominate them simultaneously in all objectives.

\subsubsection{Crossover: Generative Model based on VAE} \label{subsec:3.2.4}

\begin{figure}[H]    
  \centering         
  \includegraphics[width=\textwidth]{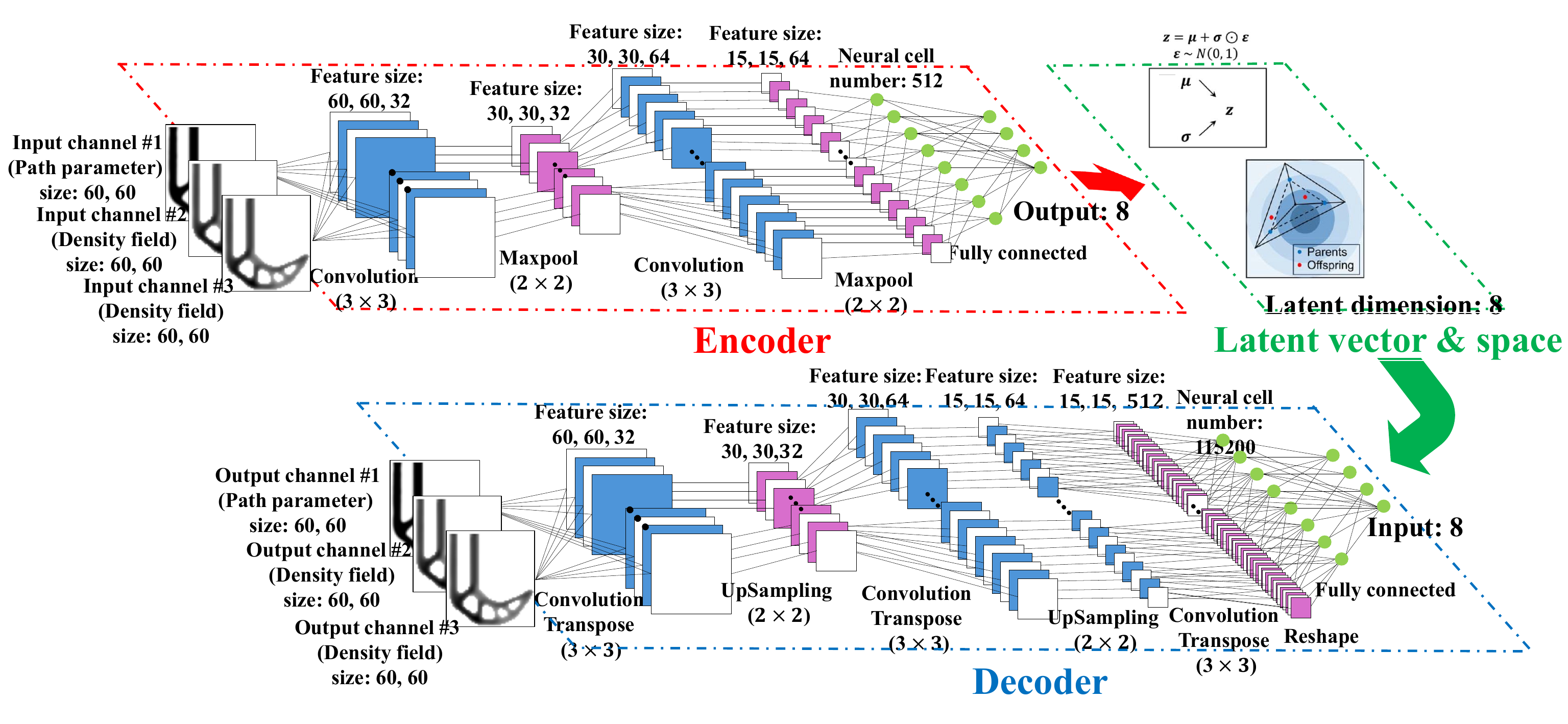}
  \caption{The architecture of multi-channel VAE.}  
  \label{fig:9}
\end{figure}

The aim in this step is to generate new candidate solutions with the characteristics of the selected elite solution in the selection step by using a generative model. We use the VAE, which is one of the representative deep generative models. The VAE model consists of two neural networks: the encoder ($p(\mathbf{z}|\mathbf{x})$) and decoder ($p(\mathbf{x}|\mathbf{z})$). 

For the encoder, the input data $\mathbf{x}$ is mapped to a probabilistic distribution in the latent space, typically assumed to be a normal distribution $q(\mathbf{z}|\mathbf{x})$. The encoder outputs the mean $\mu(\mathbf{x})$ and the variance $\sigma^2(\mathbf{x})$ of the latent variable. The objective of the encoder is to learn the distribution parameters $\mu$ and $\sigma$ in the latent space such that, given the input $\mathbf{x}$, the latent variable $\mathbf{z}$ can accurately represent the original data. While for the decoder, the latent variable $\mathbf{z}$, sampled from the latent space, is decoded back into the original data space to generate samples similar to $\mathbf{x}$. The decoder learns a conditional distribution $p(\mathbf{x}|\mathbf{z})$, generating data $\mathbf{x}$ from the latent variable $\mathbf{z}$. 
The goal of the VAE is to make the generated data $\mathbf{x'}$ as close as possible to the original data $\mathbf{x}$, while ensuring that the distribution of the latent space $q(\mathbf{z}|\mathbf{x})$ is close to a standard normal distribution $p(\mathbf{z})$. The likelihood function for the training has an explicit representation and is approximated by the evidence lower bound (ELBO):

\begin{equation} 
\mathcal{L} = \mathbb{E}_{q(\mathbf{z}|\mathbf{x})}[\log p(\mathbf{x}|\mathbf{z})] - D_\text{KL}[q(\mathbf{z}|\mathbf{x}) \parallel p(\mathbf{z})]
\label{eq:36}
\end{equation} 
 
$\mathbb{E}_{q_{\phi}(\mathbf{z}|\mathbf{x})}[\log p_{\theta}(\mathbf{x}|\mathbf{z})]$ represents the reconstruction loss, which is intuitively expected to be large, but encourages the decoder to accurately reconstruct the input data from the latent variables. From a computational point of view, this often translates to minimizing the difference between the original data and its reconstruction. $D_\text{KL}$ represents the Kullback-Leibler divergence, which is the output of the encoder. In the context of the VAE, the $D_\text{KL}$ acts as a regularizer that enforces the distribution of the latent variables to be as close as possible to the normal distribution. It also allows for the generation of new data points by sampling from the prior distribution and decoding these samples.

By balancing the reconstruction loss and the regularization term, the ELBO facilitates both effective data reconstruction and a structured distribution (often Gaussian) over the latent space. The optimization of the encoder and decoder parameters, therefore, involves minimizing the negative ELBO (hence maximizing the likelihood):

\begin{equation} 
\min (-\mathcal{L}) = - \mathbb{E}_{q(\mathbf{z}|\mathbf{x})} [\log p(\mathbf{x}|\mathbf{z})] + \beta_\text{KL} D_\text{KL} [q(\mathbf{z}|\mathbf{x}) \parallel p(\mathbf{z})]
\label{eq:37}
\end{equation} 

where the $\beta_\text{KL}$ is the weight factor to $D_\text{KL}$. This optimization ensures effective learning of the encoder and decoder parameters, leading to a VAE that proficiently captures the complexities of the data in its latent representation.

The neural network architecture of the VAE is depicted in Figure~\ref{fig:9}, we apply several convolutional layers activated by rectified linear units (ReLU) and max-pooling layers to extract features from the unit cells and recognize key patterns. Fully connected linear neural networks are used to compute the mean \(\mu\), standard deviation \(\sigma_{var}\), and latent vector \(\mathbf{z}\). The encoder transforms the input 3-channel data into a flattened dense layer, preparing it for mapping to the latent space. The decoder, which almost mirrors the encoder, utilizes convolution transpose layers with ReLU activations and upsampling layers, culminating in a convolution transpose layer with a sigmoid activation function (Sig). This structure ensures a gradual reconstruction from the latent representation.

The parameters for the convolution(ReLU), convolution transepose(ReLU), and convolution transepose(Sig) layers depend on the input image resolution and latent dimension. The latent vector’s dimension matches $\mu$ and $\sigma$, requiring the final dense layer output in the encoder to be the same size of the latent vector. During training, the decoder takes the sampled latent vector $\mathbf{z}$ (incorporating noise from $\sigma$) as input. For inference, only $\mu$ is typically used as the latent vector $\mathbf{z}$. The detail parameters for the architecture of a 3 channel each with \(60 \times 60\) pixal case is shown in Figure~\ref{fig:9}. 

It is worth noting that in the proposed framework, the deep generative model needs to be retrained in each iteration, but the size of the training set is relatively small. Therefore, we did not use a very complex network structure, sacrificing model generalization to gain computational efficiency. Finally, through a specially designed latent space sampling method \cite{yaji2024latent}, we can efficiently sample new offspring structures with parent characteristics during each crossover.

\subsubsection{Optimization performance evaluation and convergence criterion} \label{subsec:3.2.5}
In this framework, the hypervolume metric \cite{shang2020survey} is used to assess the optimization objective. It measures the volume of the objective space covered by the set of solutions, thereby effectively assessing the proximity of the Pareto front and the diversity of the solution set. Generally, the larger the hypervolume covered by the solution set, the closer these solutions are to the global optimal, indicating better optimization performance. Besides, a larger hypervolume value also indicates a broader coverage of the objective space, reflecting greater diversity in the solution set. The optimization process terminates when the hypervolume value changes by less than $0.1\%$ over 5 successive iterations while all constraints are satisfied, or after a maximum of 250 iterations.

\section{Numerical examples and discussion} \label{sec:4}

In this numerical example section, two 2D mechanical benchmarks—the L-bracket beam and the symmetric tension beam—are adopted to validate the proposed framework. We first present the design of a typical L-bracket beam and then investigate several numerical aspects, including mesh independence analysis and the effect of different population sizes. Subsequently, we systematically analyze and discuss the proposed method through the design case of a symmetric tension beam, covering the validation of the low-fidelity optimization and the MFTD optimization process. 

In this work, we use MATLAB 2022b to implement the low-fidelity optimization process and most parts of the MFTD optimization. For the remaining tasks, such as achieving CAE model realization and FEA simulation, we utilize COMSOL Multiphysics 6.2, while the implementation of the VAE is based on TensorFlow 3.6.

\subsection{L-bracket case} \label{subsec:4.1}

\begin{figure}[H]
  \centering
  \centering
  \includegraphics[width=\textwidth]{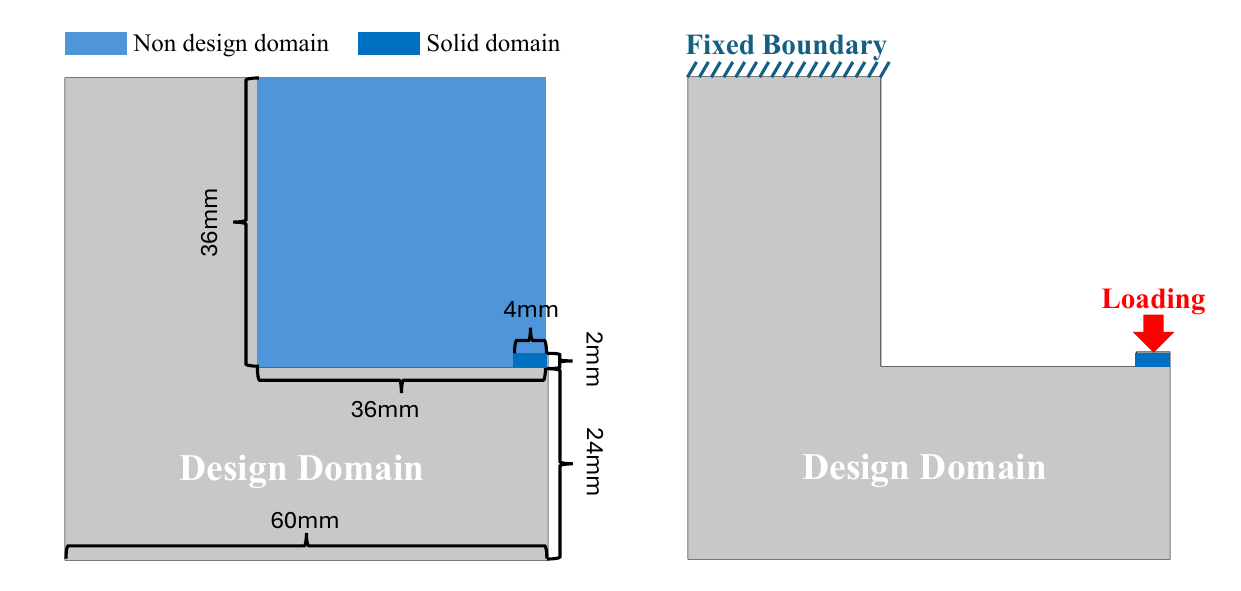} 
  \caption{The design space dimension (left) and the boundary condition (right) for the L-bracket beam case.}
  \label{fig:BC4L-bracket}
\end{figure}

The first 2D example is the design of an L-bracket part, as depicted in Figure~\ref{fig:BC4L-bracket}. To better match the data structure required for the deep model in the subsequent MFTD optimization, a square analysis domain is adopted here. The design domain and non-design domain are defined to distinguish the optimization variables. The detailed dimensions of the analysis domain are shown on the left side of Figure~\ref{fig:BC4L-bracket}. The design domain is typically defined by its ``L” shape geometry, consisting of two perpendicular arms. The top side of the vertical arm is fully fixed, and a downward vertical force 1N is applied at the right-top point of the horizontal arm (the right side of Figure~\ref{fig:BC4L-bracket}). 

\subsubsection{Parameters setting for low/high-fidelity scenario} \label{subsec:4.1.1}
A mesh consisting of \(60 \times 60\) square elements with dimensions of \(1\,\text{mm} \times 1\,\text{mm}\) is used to discretize the analysis domain for the low-fidelity optimization problem. A total of 100 low-fidelity data points are obtained based on the volume fraction \(V_\text{b}\) for the base material, ranging from 0.35 to 0.7. The detailed parameter settings for the low-fidelity optimization, high-fidelity simulation, and MFTD are shown in Table~\ref{table:TableLFHLGA}.

\begin{table}[H]
  \centering
  \caption{The table for some important optimization parameters}
  \begin{tabular}{l c c} %
    \toprule
    \textbf{Parameter} & \textbf{Symbol} & \textbf{Value} \\
    \midrule
    Filter radius for low-fidelity optimization & $R_\text{s}$ & 3 \\
    Projection Sharpness factor for low-fidelity optimization & $\beta_\text{s}$ & 32 \\
    Projection threshold value for low-fidelity optimization & $\delta_\text{s}$ & 0.5 \\
    Mass ratio for infilled solid & $V_\text{s}$ & 0.5 \\
    Penalization factor for field $\bm{\rho}$ & $p_\rho$ & 3 \\
    Penalization factor for field $\bm{\xi}$ & $p_\xi$ & 3 \\
    P-norm aggregation sharpness factor & $P$ & 12 \\
    The element dimension for each low-fidelity data & $N_\text{c}$ & 3600 \\
    Upper and downward limit of volume fraction for base material & $V_\text{b}$ & 0.35, 0.70 \\
    The morphological kernel size & $N_{\text{dilate}}$ & 4 \\
    The periodic factor of the wave function & $d$ & 8 \\
    The number of decision variables & $m$ & 100 \\
    
    The refinement factor & $\lambda$ & 9 \\
    
    \bottomrule
  \end{tabular}
  \label{table:TableLFHLGA}
\end{table}

The material parameters for the relevant materials used are listed in Table~\ref{table:TableMat}. It is worth noting that in the high-fidelity simulation, the material properties are consistent with those of the solid material used in the low-fidelity simulation.

\begin{table}[H]
  \centering
  \caption{The table for some important material properties parameters}
  \begin{tabular}{l c c} %
    \toprule
    \textbf{Parameter} & \textbf{Symbol} & \textbf{Value} \\
    \midrule
    Elastic modulus for isotropic material & $E_{\text{iso}}$ & 1 \\
    Possion ratio for isotropic materials & $v_{\text{iso}}$ & 0.3 \\
    Yield strength for isotropic material & $\sigma_\text{S}$ & 2 \\
    x/y-direction Elastic modulus for orthotropic material & $E_{11}$, $E_{22}$ & 0.3 \\
    xy-direction Elastic modulus for orthotropic material & $E_{12}$ & 0.2 \\
    xy-Possion ratio for orthotropic materials & $v_{12}$ & 0.24 \\
    x/y-direction Yield strength for isotropic material & $\sigma_\text{X}$, $\sigma_\text{Y}$ & 1 \\
    xy-direction Yield strength for isotropic material & $\sigma_\text{XY}$ & 0.5 \\
    \bottomrule
  \end{tabular}
  \label{table:TableMat}
\end{table}

\subsubsection{Mesh independence analysis} \label{subsec:4.1.2}

Mesh independence analysis is critical for stress-related problems, and many believe that mesh independence is primarily influenced by both the nesh quality and geometric shape. Therefore, we conducted a mesh independence analysis on a representative structure from our optimization case. In this work, two factors can influence the final mesh quality.

The first is the complexity of the knit curve in the generated \texttt{.dxf} file (i.e., the number of polylines). Knit curves with different levels of precision affect the accuracy of the corresponding geometry and, consequently, affect the local stress response. The second factor is the element size and meshing strategy used during adaptive mesh generation. Properly distributed element sizes improve computational accuracy, particularly in regions of stress concentration, allowing for more precise capture of stress gradients. 

Accordingly, we first categorized the derived mesh based on element size into four different modes: coarse, medium, fine, and extreme fine. The specific details are listed in Table~\ref{table:TableMeshSize}. Then, for controlling the complexity of the knit curve in the generated \texttt{.dxf} file, we utilize the Douglas-Peucker algorithm to simplify the interpolation points of the generated polylines. The primary goal of this algorithm is to reduce the number of points while preserving the overall shape and accuracy of the curves. The degree of simplification is controlled by a specified tolerance. In Table~\ref{table:TableSimplify}, we present the number of control points in the generated polylines under different tolerance values for a given geometry, as well as the corresponding number of elements in the mesh generated under the ``Fine” mesh size.

\begin{table}[H]
  \centering
  \caption{The table for the parameters of different mesh size}
  \begin{tabular}{l c c c c} %
    \toprule
    \textbf{Size Name} & \textbf{Max Size (mm)} & \textbf{Mini Size (mm)} & \textbf{Growth Rate}& \textbf{Curvature Factor}\\
    \midrule
    Coarse & 71.7 & 1.433 & 1.4 & 0.40\\
    Medium & 38.0 & 0.214 & 1.3 & 0.30\\
    Fine & 14.3 & 0.054 & 1.2 & 0.25\\
    Extreme Fine & 7.71 & 0.014 & 1.1 & 0.20\\
    \bottomrule
  \end{tabular}
  \label{table:TableMeshSize}
\end{table}

\begin{table}[H]
  \centering
  \caption{The table for the parameters with different tolerance}
  \begin{tabular}{l c c} %
    \toprule
    \textbf{Tolerance value} & \textbf{Element number} & \textbf{Number of knit} \\
    \midrule
    $10^{-1}$ & 18669 & 3091 \\
    $10^{-2}$ & 50301 & 6619 \\
    $10^{-3}$ & 62180 & 7688 \\
    $10^{-4}$ & 106351 & 11663 \\
    $10^{-5}$ & 139023 & 14809 \\
    $10^{-6}$ & 154325 & 16247 \\
    \bottomrule
  \end{tabular}
  \label{table:TableSimplify}
\end{table}

\begin{figure}[H]
  \centering
  \begin{subfigure}[b]{0.9\textwidth}
    \centering
    \includegraphics[width=\textwidth]{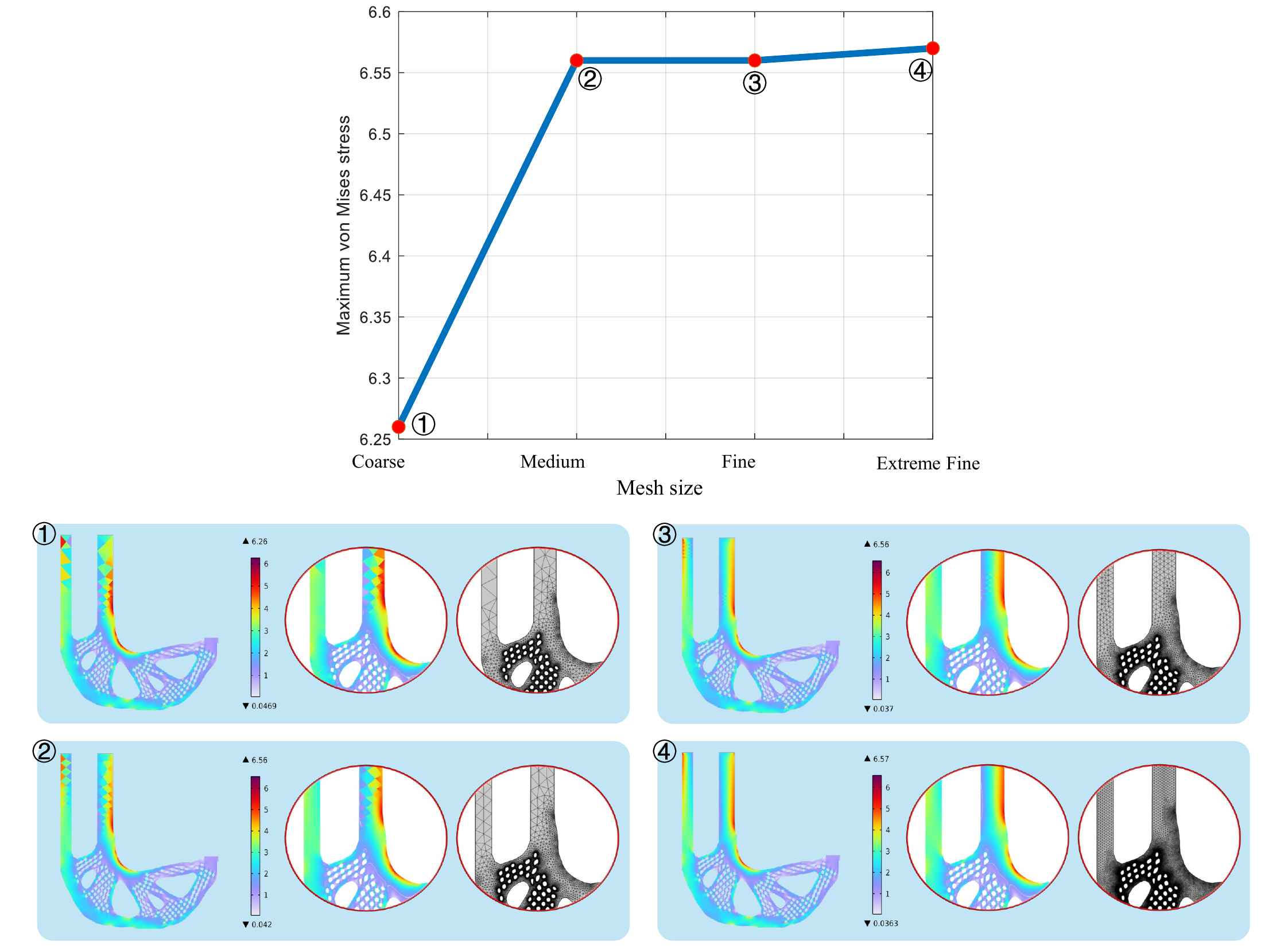} 
    \caption{}
  \end{subfigure}

  \vspace{0.5cm} 

  \begin{subfigure}[b]{0.9\textwidth}
    \centering
    \includegraphics[width=\textwidth]{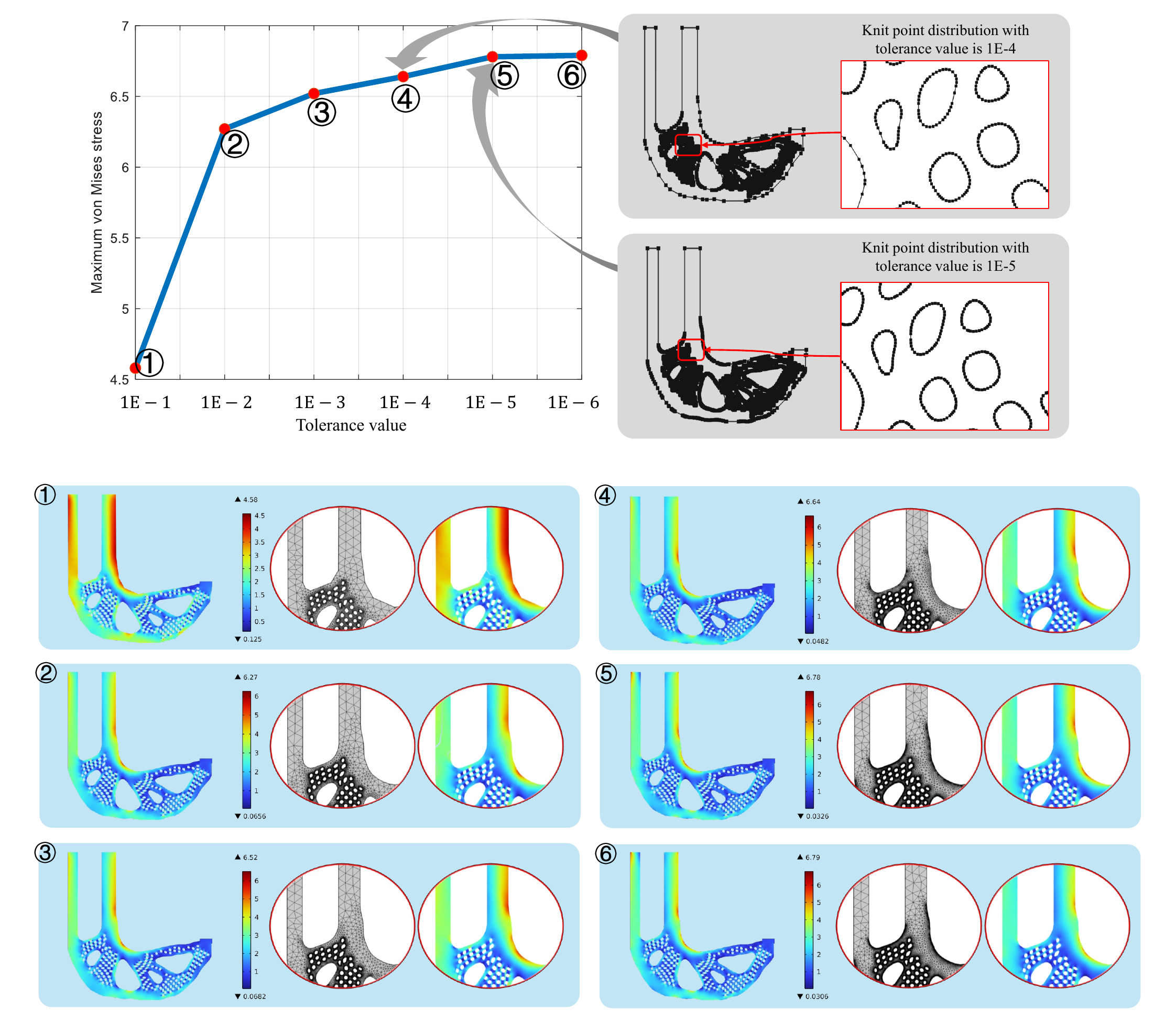} 
    \caption{}
  \end{subfigure}

  \caption{The mesh independence analysis results (a) The mesh independence analysis with respect to the element size. (b) The mesh independence analysis with respect to the tolerance value for Douglas-Peucker algorithm.}
  \label{fig:MeshIndependenceAnalysis}
\end{figure}

We first performed a mesh independence verification for different element sizes, with the results shown in the Figure~\ref{fig:MeshIndependenceAnalysis}(a). Here, we used the geometry corresponding to a tolerance value of \(10^{-4}\) as the baseline for meshing. Observations indicate that as the mesh refinement increases, the maximum von Mises stress values gradually stabilize, with minimal changes from the ``Fine” level onward. This confirms the mesh independence. Using ``Fine” or ``Extreme Fine” mesh ensures result accuracy. 

Next, we investigate the effect of the complexity of the knit curve on the maximum stress. 6 different tolerance values are used, and all cases are meshed using the ``Fine” mesh size. The specific results are shown in Figure~\ref{fig:MeshIndependenceAnalysis}(b). It can be observed that when the tolerance value is relatively large, the corresponding knit curve is significantly simplified, making it difficult to capture some detailed geometric features, and the number of generated mesh elements is relatively low. Conversely, as the tolerance value decreases, more accurate geometries are represented, which also leads to an increase in the number of mesh elements. Additionally, observations show that as the tolerance value decreases, the maximum von Mises stress values gradually stabilize, with minimal changes from \(10^{-5}\) downward. This also confirms the mesh independence.

Finally, to ensure geometric accuracy and mesh independence while minimizing computational cost, we set the mesh size to ``Fine” and the tolerance in the Douglas-Peucker algorithm to \(10^{-5}\) for all subsequent cases.

\subsubsection{Results obtained from the proposed method} \label{subsec:4.1.3}

\begin{figure}[H]
  \centering
   \begin{subfigure}[b]{1\textwidth}
    \centering
    \includegraphics[width=\textwidth]{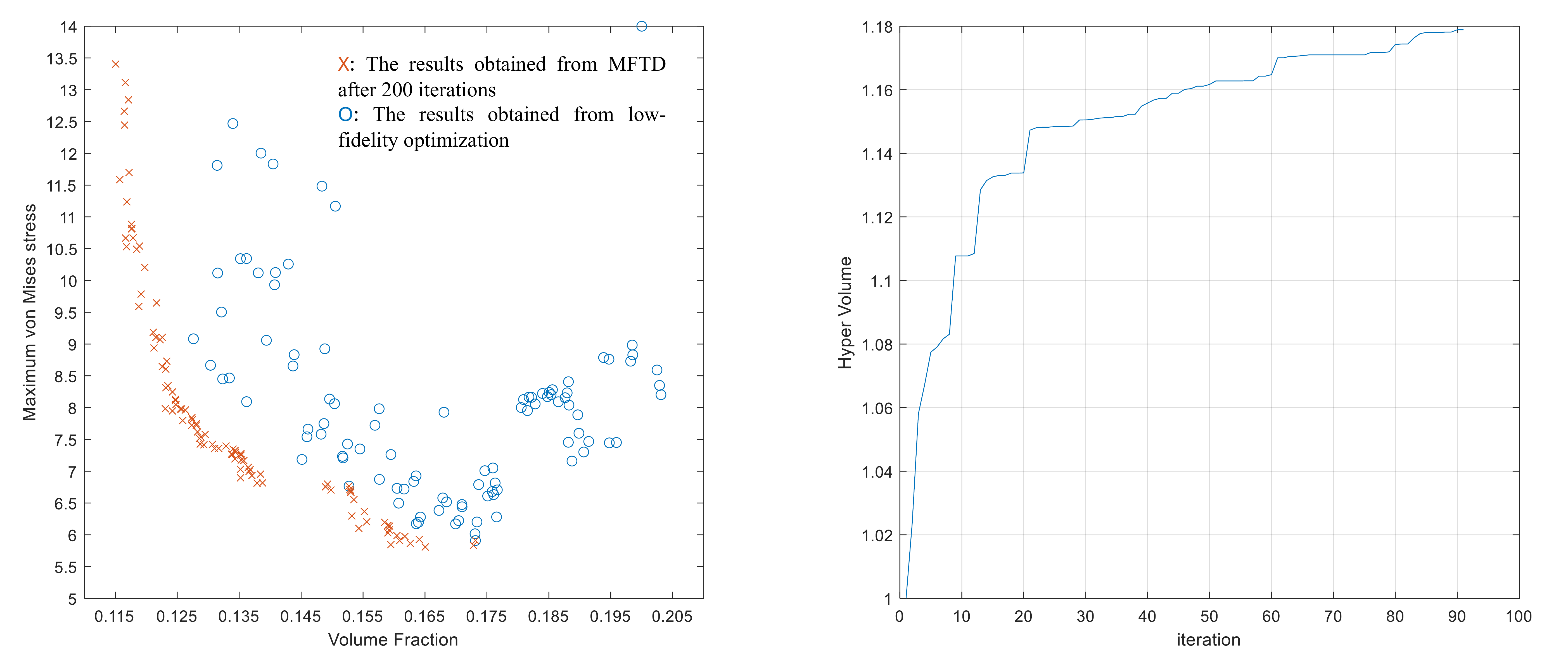} 
    \caption{}
  \end{subfigure}

  \vspace{0.5cm} 

  \begin{subfigure}[b]{1\textwidth}
    \centering
    \includegraphics[width=\textwidth]{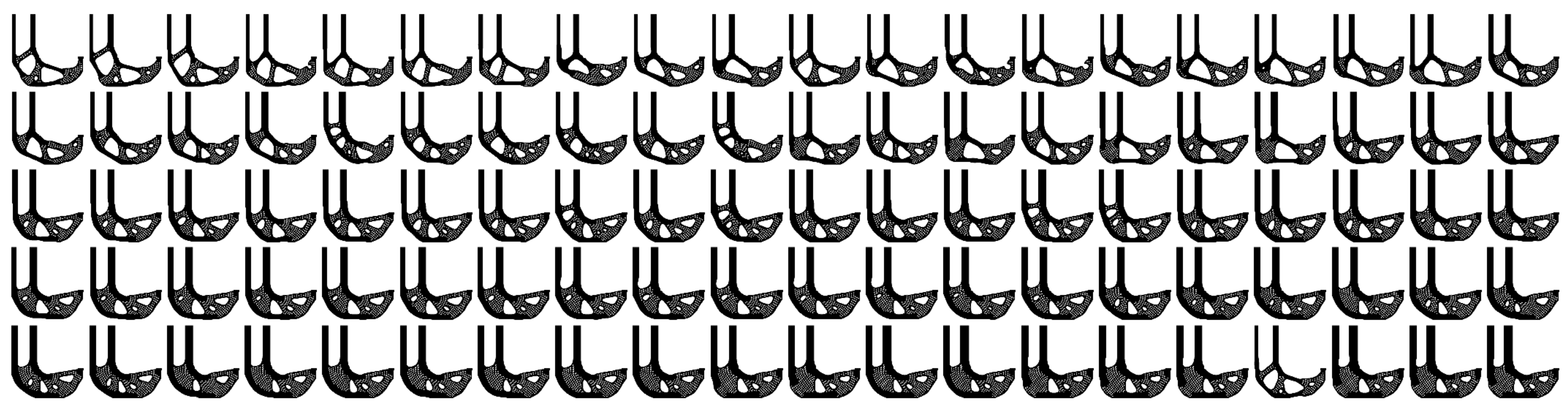} 
    \caption{}
  \end{subfigure}

  \vspace{0.5cm} 

  \begin{subfigure}[b]{1\textwidth}
    \centering
    \includegraphics[width=\textwidth]{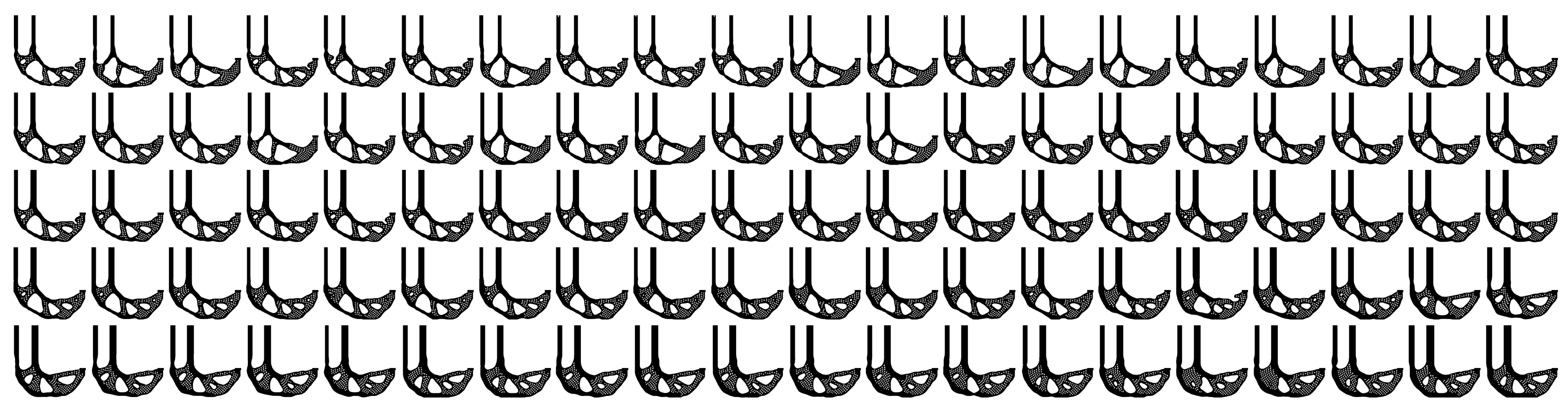} 
    \caption{}
  \end{subfigure}

  \caption{(a) Left: the optimized results between $1^\text{st}$ iteration and $90^\text{th}$ iteration; Right: the hypervolume value convergence histroy; (b) the initial 100 high-fidelity results; (c) the final 100 optimized high-fidelity results.}
  \label{fig:100results4L-bracket}
\end{figure}

Following the high-fidelity realization workflow shown in Figure~\ref{fig:MFTDflowchat}, 100 initial high-fidelity results are obtained (see Figure~\ref{fig:100results4L-bracket}(b)). We ultimately derive the elite results displayed in Figure~\ref{fig:100results4L-bracket}(c). The optimization is stopped at the 90\textsuperscript{th} iteration. Figure~\ref{fig:100results4L-bracket}(a) demonstrates the improvement in the performance of the elite material distributions. The circles represent the initial results, while the red crosses indicate the final elite solutions. It can be observed that the elite results are converging towards minimization in the objective function space, confirming that the performance of the elite results improves throughout the iterations. As shown in Figure~\ref{fig:100results4L-bracket}(a) left, the hypervolume indicator reaches approximately 1.18 at the end of the optimization, indicating that the overall performance of the elite results has improved significantly.

\begin{figure}[H]    
  \centering         
  \includegraphics[width=\textwidth]{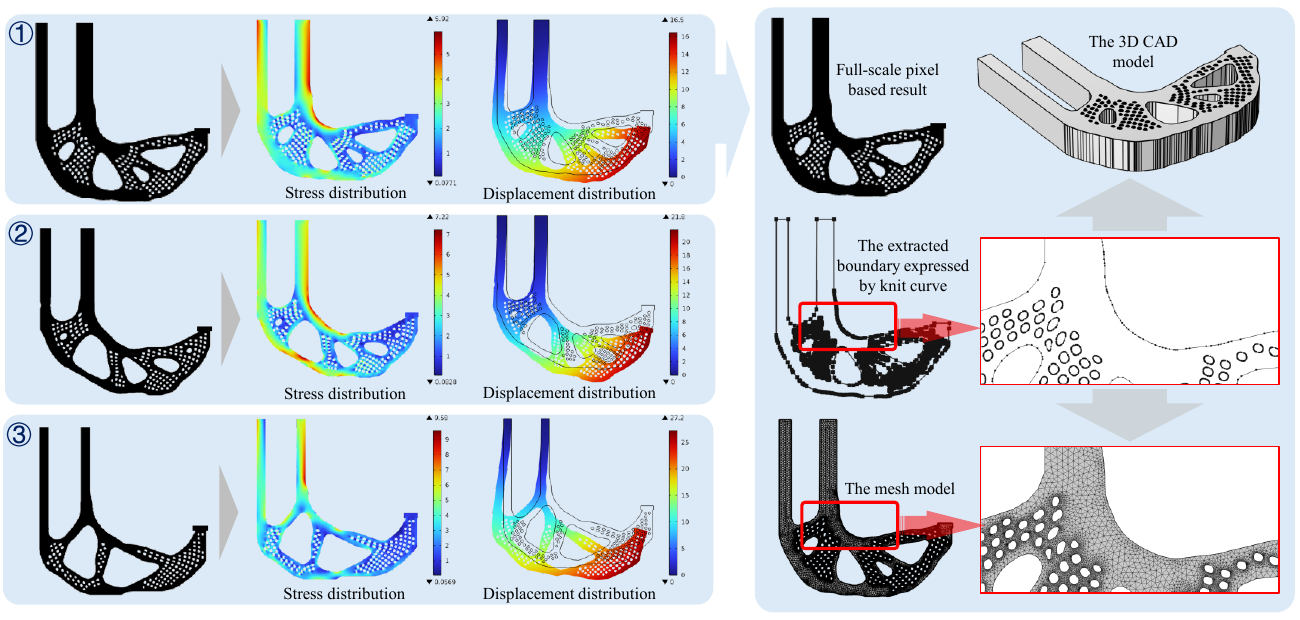}
  \caption{The optimized results from the MFTD framework for L-bracket beam case.}  
  \label{fig:Results4Lbracket100}
\end{figure}

To better evaluate the optimized results, we select the results with mass ratio \(G_\text{vf}\) of 15.6\%, 13.4\%, and 11.8\% as examples (see Figure~\ref{fig:Results4Lbracket100}(a)), and then analyze it by focusing on the structural deformation and von Mises stress distribution. It could be observed that both the optimized L-bracket beams have complex geometry features. The maximum displacement occurs near the right end of the horizontal arm (the three far-right column shown in Figure~\ref{fig:Results4Lbracket100}), where the applied load is likely highest, reaching values of around 16.5mm, 21.8mm, and 27.2mm respectively, as indicated by the red color on the scale. This suggests significant bending and deformation in the unsupported region due to the applied force. 

Unlike the similarity observed in the deformation distribution, the stress distribution of the three results exhibits significant diversity. For the result with \(G_\text{vf} = 15.55\%\), the high stress regions are primarily concentrated at the connection corners of the supports and near the base of the loaded column. These areas are likely the main regions of stress concentration or bending within the structure and represent its weak points; hence, they are filled with solid material. For the result with \(G_\text{vf} = 13.39\%\), the inner and outer sides of the L-shaped corners exhibit the highest stress concentration. This is because these locations experience significant bending and shear stress under loading. As a typical stress concentration point, these areas are also filled with solid material. For the result with \(G_\text{vf} = 11.85\%\), the junction between the vertical column base and the horizontal beam, as well as the inner surface of the vertical column on the inside, exhibits higher stress distribution. Compared to the previous two results, this structure has the thinnest vertical column, which bears the majority of the load, resulting in significant shear and bending stress at the connection points. Overall, it can be observed that in all three structures, the stress distribution around the porous infill regions is relatively uniform, indicating that the hole design is reasonable and does not cause significant stress concentration.

\subsubsection{Comparison amongst different decision variable number} \label{subsec:4.1.4}

Within MFTD, the population size greatly influences the final optimization results. In this work, the population size will affect the quality of VAE training and, consequently, the quality of the generated results. Therefore, in this section, we investigate the impact of population size by varying the number of decision variables in the case shown in Figure~\ref{fig:BC4L-bracket}.

\begin{figure}[H]    
  \centering
   \begin{subfigure}[b]{1\textwidth}
    \centering
    \includegraphics[width=\textwidth]{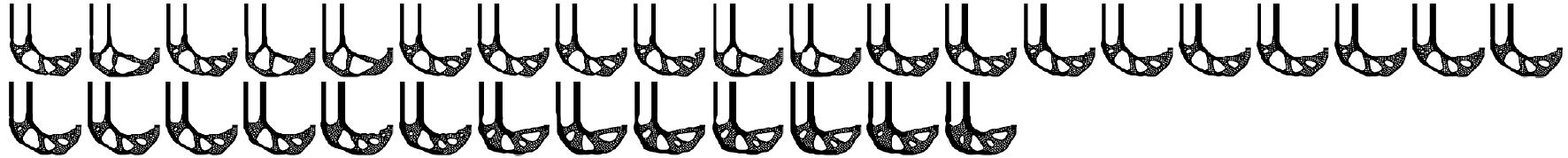} 
    \caption{}
  \end{subfigure}

  \vspace{0.5cm} 

  \begin{subfigure}[b]{1\textwidth}
    \centering
    \includegraphics[width=\textwidth]{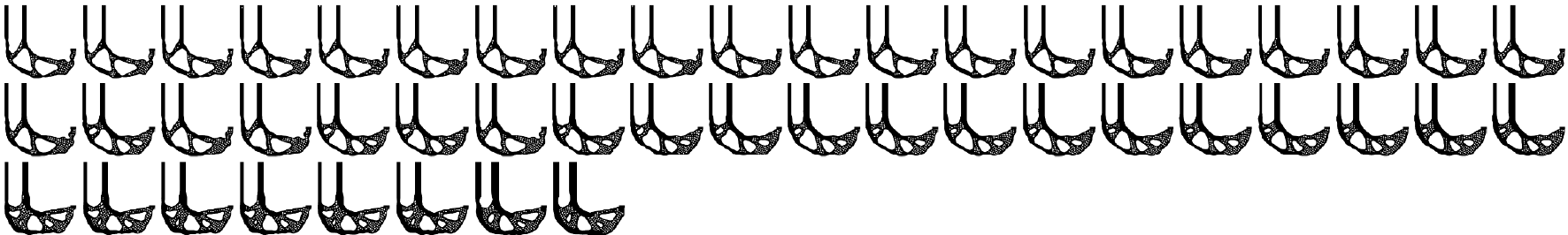} 
    \caption{}
  \end{subfigure}

  \vspace{0.5cm} 

  \begin{subfigure}[b]{1\textwidth}
    \centering
    \includegraphics[width=\textwidth]{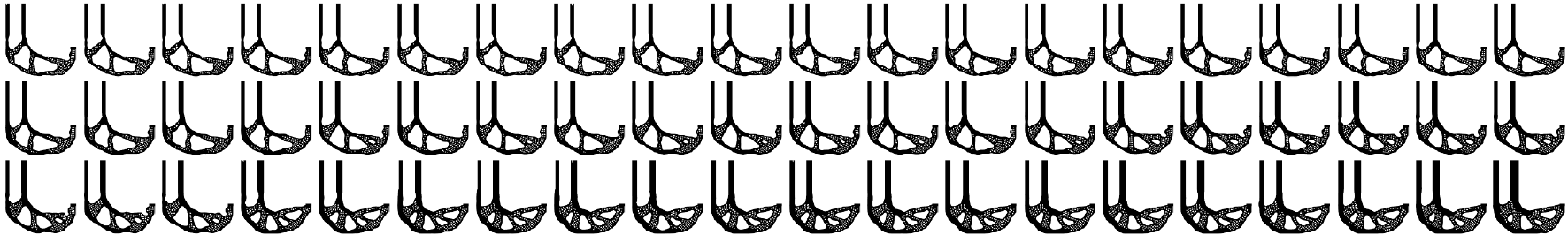} 
    \caption{}
  \end{subfigure}

  \caption{The final obtained Rank-1 results (a) The 100-results; (b) The 150-results; (c) The 200-results.}  
  \label{fig:Rank1-4Lbracket}
\end{figure}

We additionally performed MFTD optimization using two additional population sizes of 150 and 200, with identical optimization and material parameters. High-fidelity analysis was performed on this series of results, and the Rank-1 results were ultimately selected. For the results with a population size of 100 (100-results), the number of Rank-1 results obtained was 33. In contrast, the results with population sizes of 150 (150-results) and 200 (200-results) yielded 48 and 60 Rank-1 results, respectively. The detailed structures are shown in Figure~\ref{fig:Rank1-4Lbracket}.

\begin{figure}[H]
  \centering
  \begin{subfigure}[b]{0.45\textwidth}
    \centering
    \includegraphics[width=\textwidth]{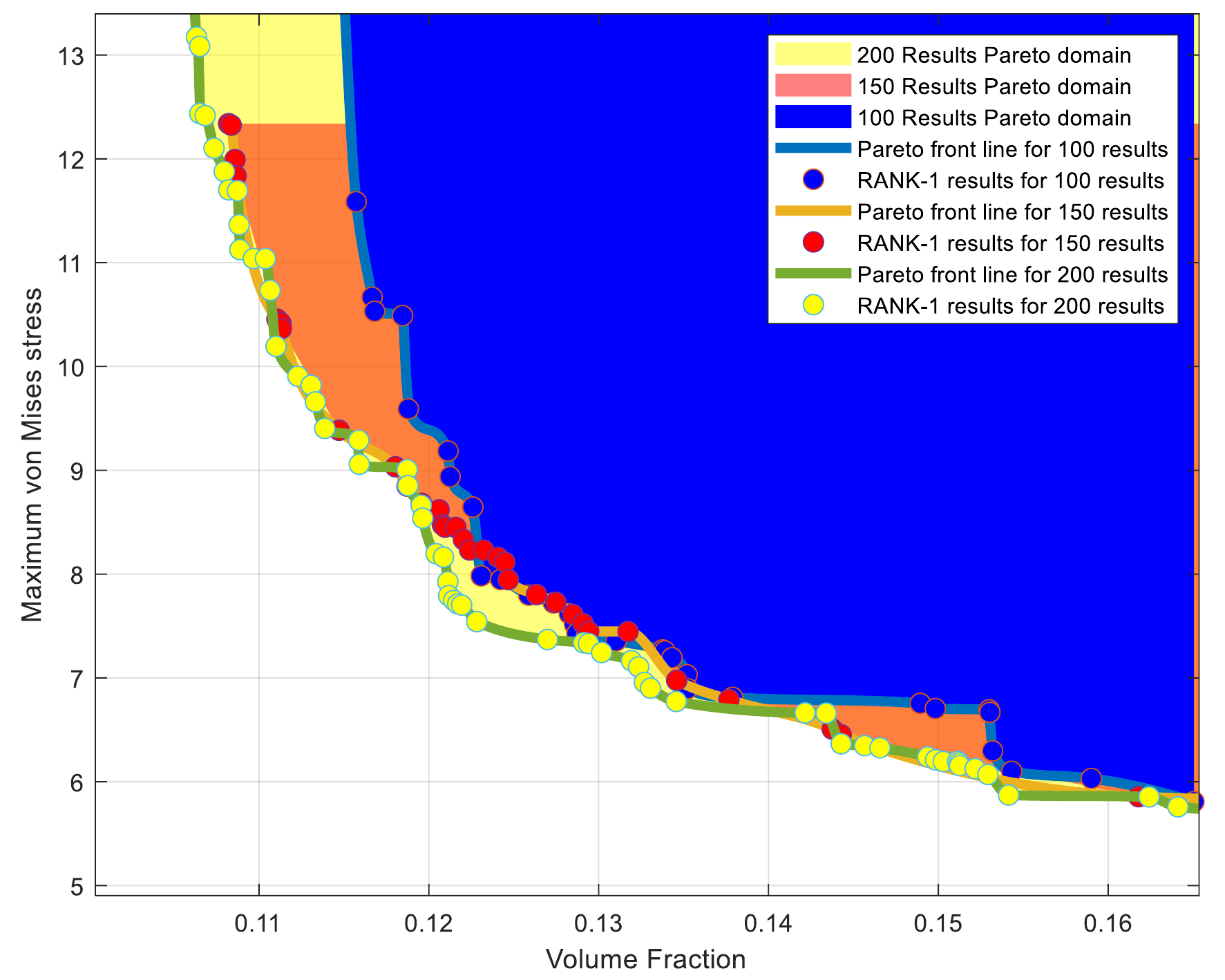} 
    \caption{}
  \end{subfigure}
  \hfill
  \begin{subfigure}[b]{0.51\textwidth}
    \centering
    \includegraphics[width=\textwidth]{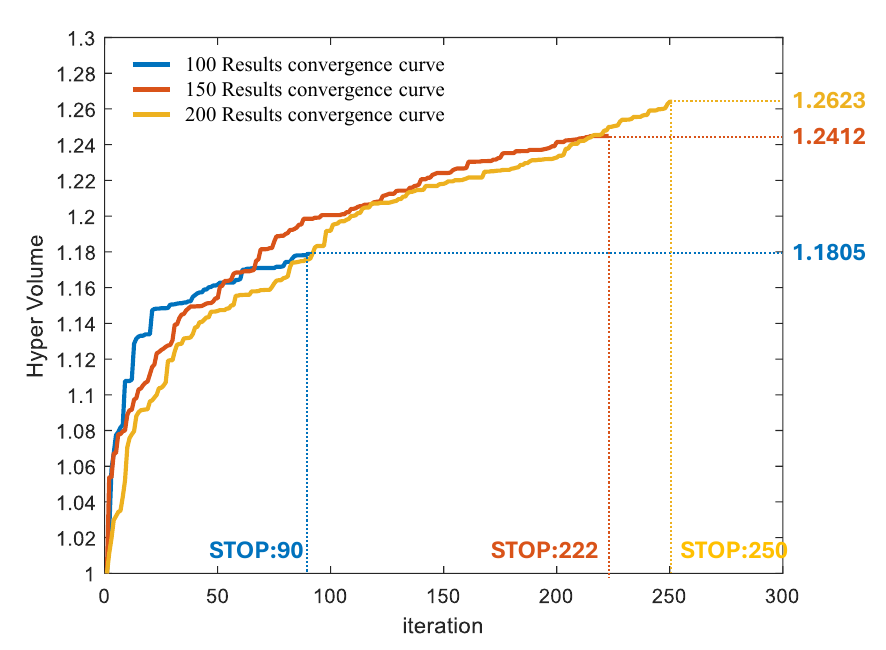} 
    \caption{}
  \end{subfigure}
  \caption{The Pareto front domain for three different optimization strategies (a); and their corresponding convergence histroy (b).}  
  \label{fig:ParetoResults4Lbracket}
\end{figure}

The Pareto front domain formed by the Rank-1 solutions for the three optimizations is illustrated in Figure~\ref{fig:ParetoResults4Lbracket}(a). Clearly, as the population size increases, the Pareto front domain hypervolume also becomes larger, indicating the potential to generate better results. However, compared to the 150-results, the improvement in the 200-results is not very significant, with noticeable enhancement only within a specific volume fraction range (12\%~13\%). The corresponding hypervolume convergence curves are shown in Figure~\ref{fig:ParetoResults4Lbracket}(b). Observations reveal that the 100-results optimization process stopped at the $90^\text{th}$ iteration, while the 150-results and 200-results continued until the $222^\text{nd}$ and $250^\text{th}$ iterations, respectively. Notably, for the 200-results optimization, the process had not yet reached the specified convergence criterion even at the maximum iteration step of 250.

\begin{figure}[H]    
  \centering         
  \includegraphics[width=\textwidth]{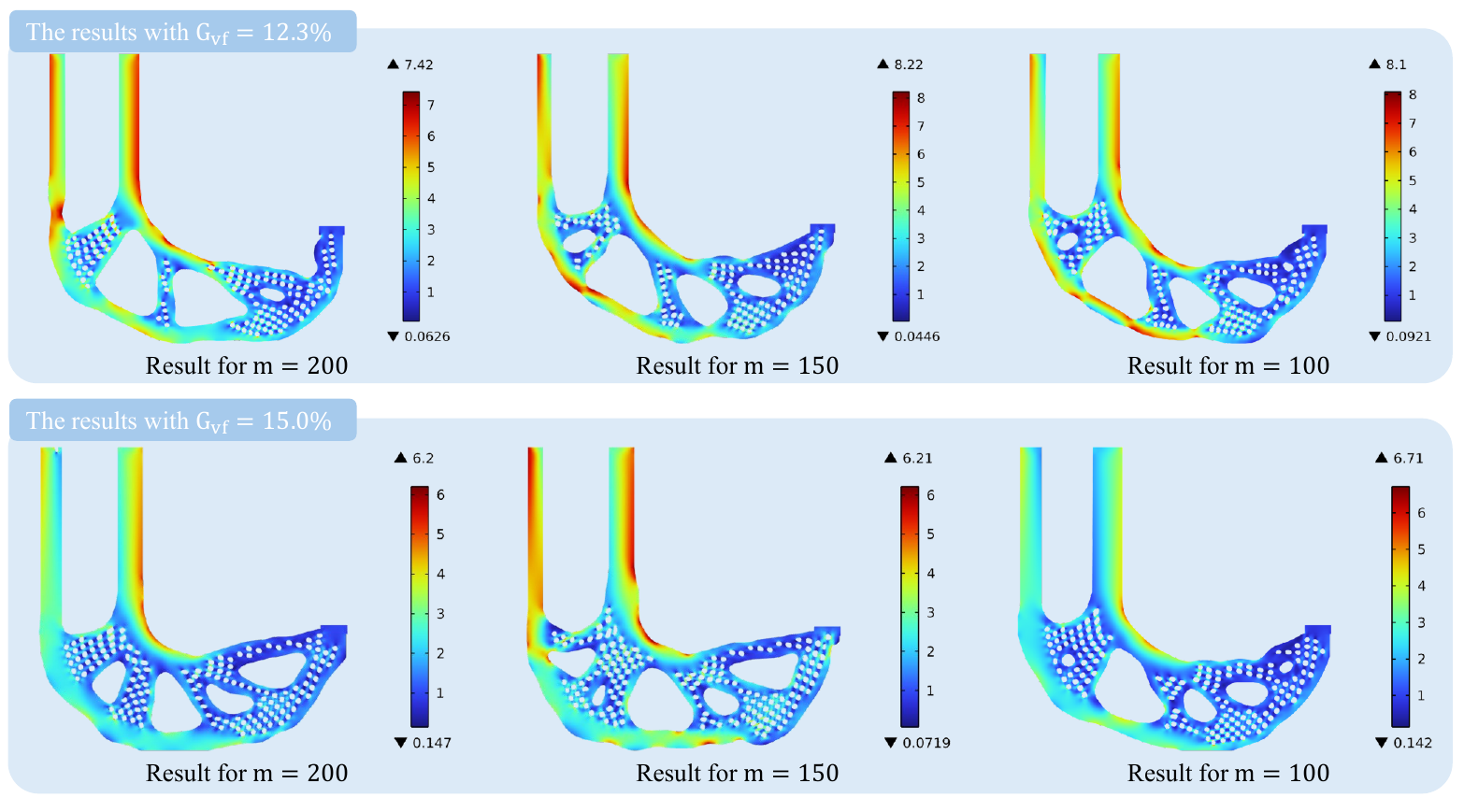}
  \caption{The selected structures obtained by three different optimization strategies with specific volume fractions. (a) \(G_\text{vf} = 12.3\%\); (b) \(G_\text{vf} = 15.0\%\).}  
  \label{fig:Structures4LF}
\end{figure}

To further facilitate a visual comparison, we selected two sets of results with volume fractions of 12.3\% and 15\% for comparison, and the corresponding results are shown in Figure~\ref{fig:Structures4LF}. As seen in Figure~\ref{fig:Structures4LF}, for both sets of results, the structures generated by the 200-results exhibit the lowest von Mises stress. Notably, in the \(G_\text{vf} = 12.3\%\) group, the 200-results produced a completely different topology compared to the other two optimizations. In the \(G_\text{vf}= 15.0\%\) group, the structures generated by the 200-results and 150-results exhibit similar performance (with nearly identical maximum stress values), but their topological configurations are distinct. 

Additionally, it is worth noting that the base material distribution for the 100-results in both the 12.3\% and 15\% groups is highly similar. In contrast, the results generated by the 200-results and 150-results in the two groups display significant differences. Since no mutation operations were incorporated into this framework, this demonstrates the impact of population size on population diversity.

\subsection{Symmetric tension beam case} \label{subsec:4.2}
\begin{figure}[H]
  \centering
  \centering
  \includegraphics[width=\textwidth]{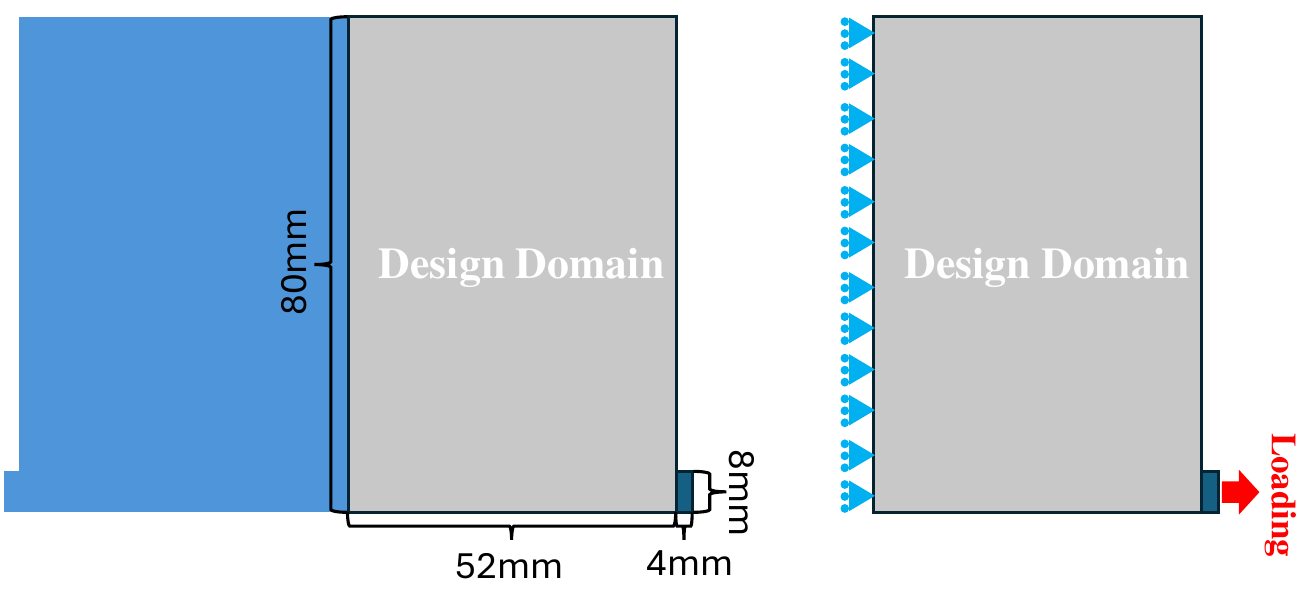} 
  \caption{The design space dimension (left) and the boundary condition (right) for the symmetric tension beam case.}
  \label{fig:BC4symbeam}
\end{figure}

The second 2D example involves designing a symmetric tension beam. In this section, we further validate the effectiveness of our method. The boundary conditions and structural dimensions are illustrated in Figure~\ref{fig:BC4symbeam}. The beam has a rectangular shape, with its width slightly greater than its length. The beam is symmetric about its vertical centerline, allowing for optimization of only half of the beam to reduce computational complexity. The degrees of freedom on the left side are fixed in the \(x\)-direction, restricting movement to the vertical direction while allowing free horizontal displacement. A horizontal outward force of 1~N is applied at the bottom-right corner.

\subsubsection{Parameters setting for low/high-fidelity scenario} \label{subsec:4.2.1}
For the low-fidelity optimization problem, a mesh consisting of \(52 \times 80\) square elements with dimensions of \(1\,\text{mm} \times 1\,\text{mm}\) is utilized to discretize the design domain. The parameter settings for low-fidelity optimization are the same as those listed in Table~\ref{table:5}. Only few parameters for the high-fidelity simulation and MFTD differ from those used in the L-bracket case, as detailed in the table below:

\begin{table}[H]
  \centering
  \caption{The table for some new optimization parameters for the symmetric tension beam case}
  \begin{tabular}{l c c} %
    \toprule
    \textbf{Parameter} & \textbf{Symbol} & \textbf{Value} \\
    \midrule
    The morphological kernel size & $N_{\text{dilate}}$ & 6 \\
    The periodic factor of the wave function & $d$ & 10 \\
    The number of decision variables & $m$ & 160 \\
    Upper and downward limit of volume fraction & $V_\text{s}$ & 0.3, 0.8 \\
    The refinement factor & $\lambda$ & 16 \\
    The element dimension for each low-fidelity data & $N_\text{c}$ & 3200 \\
    \bottomrule
  \end{tabular}
  \label{table:5}
\end{table}

Other optimization parameters and relevant material properties are identical to those used in the L-bracket case.

\subsubsection{Effectiveness of low-fidelity optimization} \label{subsec:4.2.2}

Here, we first validate the effectiveness of the low-fidelity optimization proposed in Section 2.1. To this end, we compare a series of results with varying volume fractions obtained from low-fidelity optimization with those from traditional single-material, stress-based topology optimization (Stress-TO). To ensure a fair comparison, all optimization or material parameters are kept consistent. Subsequently, we apply the method proposed in Section 2.2 to transform these results and obtain the corresponding high-fidelity outcomes.

\begin{figure}[H]    
  \centering         
  \begin{subfigure}[b]{1\textwidth}
   \centering
   \includegraphics[width=\textwidth]{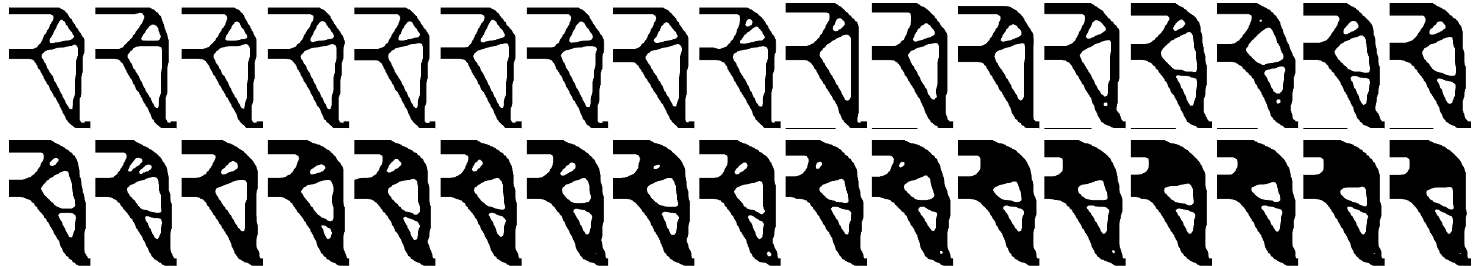} 
   \caption{}
 \end{subfigure}

 \vspace{0.5cm} 

 \begin{subfigure}[b]{1\textwidth}
   \centering
   \includegraphics[width=\textwidth]{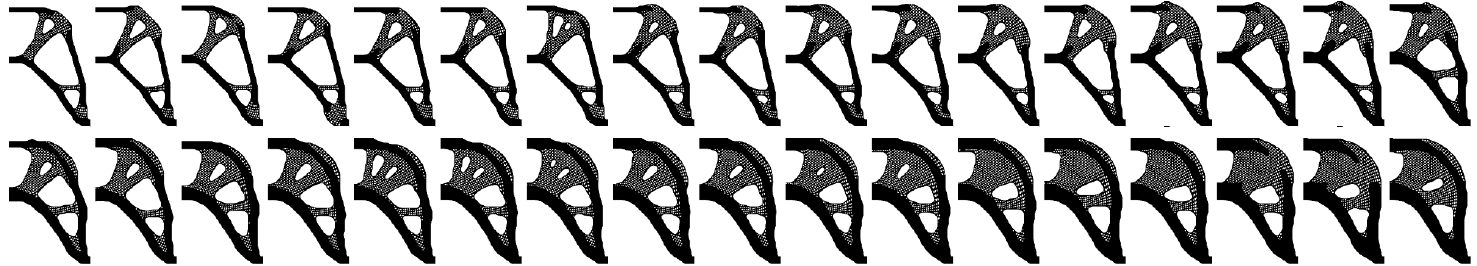} 
   \caption{}
 \end{subfigure}

  \caption{The final obtained Rank-1 results (a) Stress-TO; (b) low-fidelity optimization.}  
  \label{fig:Rank1stressLF}
\end{figure}

\begin{figure}[H]    
  \centering         
  \begin{subfigure}[b]{1\textwidth}
    \centering
    \includegraphics[width=\textwidth]{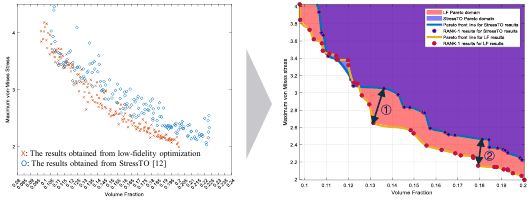} 
    \caption{}
  \end{subfigure}

  \vspace{0.5cm} 
  
  \begin{subfigure}[b]{1\textwidth}
    \centering
    \includegraphics[width=\textwidth]{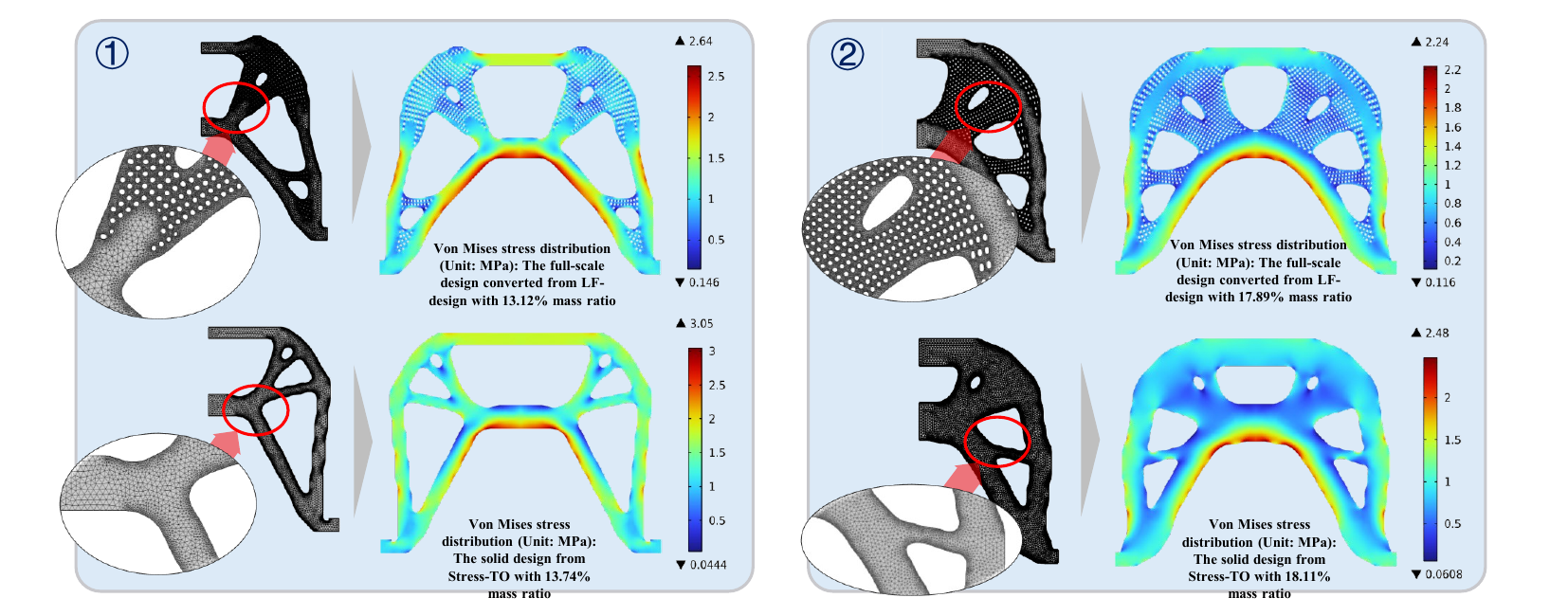} 
    \caption{}
  \end{subfigure}

  \caption{The structural performance comparison for the low-fidelity optimization and Stress-TO (a); and two sets of results with detail structure with volume fractions of 13\% and 18\% for comparison (b).}  
  \label{fig:ParetoStressLF}
\end{figure}

Within the volume fraction range of 0.3 to 0.8, we generated 160 results using both the low-fidelity optimization and Stress-TO methods. High-fidelity analysis was performed on this series of results, and the Rank-1 results were ultimately selected. Coincidentally, both the low-fidelity optimization and Stress-TO methods produced 34 Rank-1 results. The detailed structures are shown in Figure~\ref{fig:Rank1stressLF}. 
The corresponding Pareto results are shown in the left side of Figure~\ref{fig:ParetoStressLF}, and the Pareto domain formed by the Rank-1 solutions is illustrated in the right side of Figure~\ref{fig:ParetoStressLF}. The Pareto domain hypervolume for low-fidelity optimization (the area of the red region) is larger than that of Stress-TO (the area of the blue region), indicating that low-fidelity optimization utilizes the solution space more effectively and generates better results compared to Stress-TO.

To further facilitate a visual comparison, we selected two sets of results with volume fractions of 13\% and 18\% for comparison, and the corresponding results are shown in Figure~\ref{fig:ParetoStressLF} (b). As seen in the figure, for the group with a volume fraction of 13\%, the structure generated by Stress-TO has a volume fraction of 13.74\% and a maximum stress of 3.05 MPa. In contrast, the structure generated by low-fidelity optimization has a volume fraction of 13.12\% and a maximum stress of 2.64 MPa, representing a 13\% reduction in maximum stress. For the group with a volume fraction of 18\%, the structure generated by Stress-TO has a volume fraction of 18.11\% and a maximum stress of 2.48 MPa, while the structure generated by low-fidelity optimization has a volume fraction of 17.89\% and a maximum stress of 2.24 MPa, representing a 9.7\% reduction in maximum stress.

\subsubsection{Results obtained from the proposed method} \label{subsec:4.2.3}

\begin{figure}[H]    
  \centering         
  \begin{subfigure}[b]{1\textwidth}
    \centering
    \includegraphics[width=\textwidth]{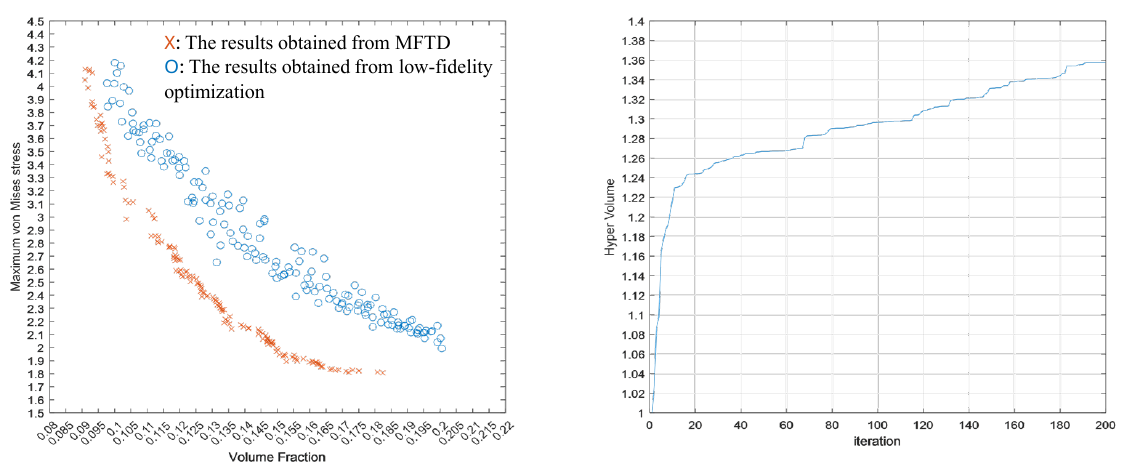} 
    \caption{}
  \end{subfigure}

  \vspace{0.5cm} 

  \begin{subfigure}[b]{1\textwidth}
    \centering
    \includegraphics[width=\textwidth]{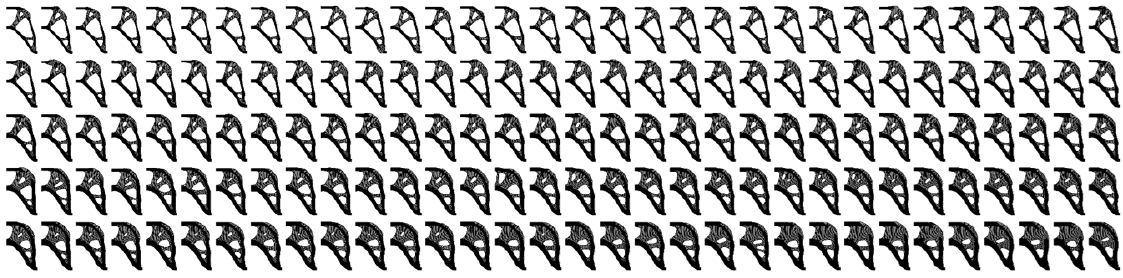} 
    \caption{}
  \end{subfigure}

  \vspace{0.5cm} 

  \begin{subfigure}[b]{1\textwidth}
    \centering
    \includegraphics[width=\textwidth]{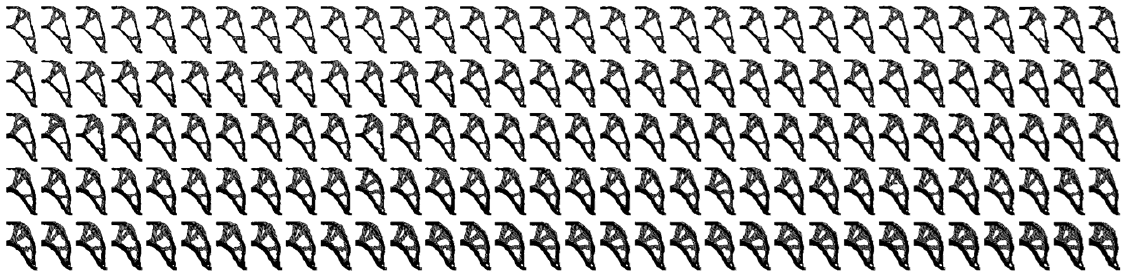} 
    \caption{}
  \end{subfigure}

  \caption{(a) Left: the optimized results between $1^\text{st}$ iteration and $200^\text{th}$ iteration; Right: the hypervolume value convergence histroy; (b) the original 160 low-fidelity optimization results; (c) the 160 MFTD results after 200 iterations.}  
  \label{fig:MFTDresults4case2}
\end{figure}

Building on the 160 initial solutions obtained using low-fidelity optimization in the previous section, we further performed MFTD optimization on them. Figure~\ref{fig:MFTDresults4case2} presents several key aspects of the optimization process and results for the topology design. Figure~\ref{fig:MFTDresults4case2}(a) illustrates the Pareto front obtained from the optimization process, comparing volume fraction and the optimization number across iterations. The embedded images display representative designs along the Pareto front, showing the evolution of material distribution through the optimization process. The convergence history of the hypervolume value is shown in Figure~\ref{fig:MFTDresults4case2}(b), and the optimization process finally stopped at the $200^\text{th}$ iteration with a increase of 35.5\%. 
Figure~\ref{fig:MFTDresults4case2}(c) displays 98 results which are finally selected as the Rank-1 results.

\begin{figure}[H]
  \centering         
  
  \begin{subfigure}[b]{1\textwidth}
    \centering
    \includegraphics[width=\textwidth]{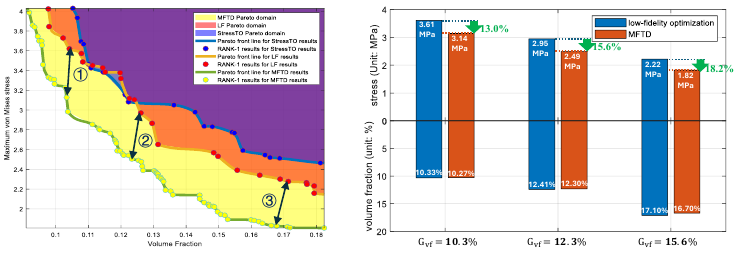} 
    \caption{}
  \end{subfigure}

  \vspace{0.5cm} 

  \begin{subfigure}[b]{0.8\textwidth}
    \centering
    \includegraphics[width=\textwidth]{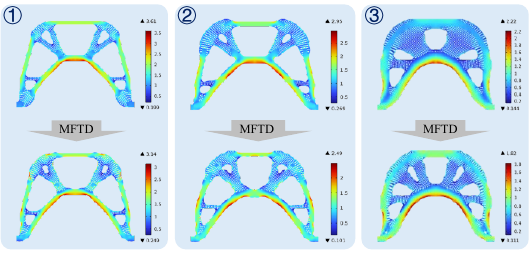} 
    \caption{}
  \end{subfigure}

  \caption{The structural performance comparison for the MFTD, low-fidelity optimization, and Stiff-TO for the symmetric tension beam case (a); and the detail structure comparison between MFTD and low-fidelity optimization (b).}  
  \label{fig:Pareto4case2}
\end{figure}

The MFTD Pareto domain formed by the Rank-1 solutions is illustrated in the Figure~\ref{fig:Pareto4case2} left part and is compared with the results shown in Figure~\ref{fig:ParetoStressLF}. Clearly, the hypervolume for MFTD (the area of the yellow region) is larger than that of Stress-TO (the area of the blue region) and low-fidelity optimization (the area of the orange region), indicating that the MFTD results are the best among the three. To further demonstrate the effectiveness of the MFTD optimization process, the results with mass ratios \(G_\text{vf}\) of 10.3\%, 12.3\%, and 17\% for both the low-fidelity optimization and MFTD results were selected as examples. The corresponding structures and their von Mises stress distributions are shown on the right side of Figure~\ref{fig:Pareto4case2}(a).

It can be observed that for the group with \(G_\text{vf} = 10.3\%\), the structure generated by low-fidelity optimization has a volume fraction of 10.33\% and a maximum stress of 3.61 MPa. In contrast, the structure generated by MFTD has a volume fraction of 10.27\% and a maximum stress of 3.14 MPa, achieving a 13\% reduction in maximum stress with a smaller volume fraction. For the group with \(G_\text{vf}= 12.3\%\), the results of low-fidelity optimization and MFTD have volume fractions of 12.3\% and 12.41\%, respectively, and maximum stresses of 2.95 MPa and 2.49 MPa. The MFTD design achieves a 0.9\% reduction in weight and a 15.6\% improvement in maximum stress. Finally, for the group with \(G_\text{vf} = 17.0\%\), the results of low-fidelity optimization and MFTD have volume fractions of 16.7\% and 17.1\%, respectively, and maximum stresses of 2.22 MPa and 1.82 MPa. The MFTD design achieves a 2.3\% reduction in weight while improving the reduction of maximum stress by 18.2\%.

We subsequently analyzed the geometric configurations of the three sets of results shown in Figure~\ref{fig:Pareto4case2}(b). First, regarding the distribution of the base material, for the groups with \(G_\text{vf} = 10.3\%\) and \(G_\text{vf} = 12.3\%\), the topology of the structures generated by low-fidelity optimization and MFTD remains unchanged, with differences only in shape. However, for the group with \(G_\text{vf} = 17.0\%\), both the topology and shape differ between the two methods. Regarding the solid-filled regions, significant differences are observed between the structures generated by low-fidelity optimization and MFTD across all three groups. A comprehensive observation indicates that, compared to the structures from low-fidelity optimization, the solid-filled regions in the MFTD structures are more continuous, better aligning with the overall load-bearing directions of the structure.

\section{Conclusion} \label{sec:5}
In this paper, we utilize a data-driven approach within the MFTD framework to achieve a hybrid solid-porous infill design that considers the stress concentration issue. Within this framework, the detailed topology of the complex hybrid solid-porous infill structure is directly controlled by three control fields. These fields are initialized through a specifically designed low-fidelity optimization process and subsequently optimized using the MFTD method. The optimized result from this approach features a detailed and accurate geometric shape, with the performance of the precise model considered throughout the optimization process. Numerical results indicate that the developed method can robustly generate innovative designs across various cases, ensuring that solid infill materials are distributed at critical structural positions.

Based on the discussions and analyses presented in the previous sections, the following conclusions can be drawn:

1. By carefully controlling the geometric mapping accuracy and mesh size, we ensure stress analysis mesh convergence for the designed geometric patterns. This guarantees the reliability and accuracy of subsequent optimizations.

2. The side-by-side mapping between control variables and geometric structures allows us to analyze and obtain precise CAD models during the optimization process. These models can be directly used for prototype fabrication, effectively reducing the gap between design and manufacturing and achieving a design-for-manufacturing approach.

3. The population size significantly impacts the optimization results. Increasing the population size as much as possible, within the constraints of computational cost, enhances the quality of the results.

4. Due to the nonlinearity of stress analysis and its localized response characteristics, it is challenging to achieve a Rank-1 result count close to the total population size within a finite number of iterations.

5. Although the proposed low-fidelity optimization incorporates several simplifications, it still demonstrates performance advantages compared to traditional single-material stress optimization. When used as an initial solution for the MFTD optimization process, these advantages are further amplified, highlighting the value of the proposed framework.

There are several limitations and areas for further exploration:

1. Currently, the developed method is focused solely on 2D problems, and future work will extend it to 3D practical engineering applications.

2. Additionally, this study did not consider buckling performance in the optimization process, which will be addressed in future work.

3. We have not yet conducted experimental validation for the fabricated prototypes, and this will be completed in subsequent studies.








\section*{Acknowledgements}
This work was supported by JSPS KAKENHI, Grant Number 23H03799. 

\bibliographystyle{unsrt} 
\bibliography{Reference} 

\end{document}